# Graphon Mean-Field Logit Dynamic: Derivation, Computation, and Applications

Hidekazu Yoshioka[1, *]

[1] Graduate School of Advanced Science and Technology, Japan Advanced Institute of Science and Technology, 1-1 Asahidai, Nomi, Ishikawa, Japan

[*] Corresponding author: yoshih@jaist.ac.jp, ORCID: 0000-0002-5293-3246

**Abstract**

We present a graphon mean-field logit dynamic, a stationary mean-field game based on logit interactions. This dynamic emerges from a stochastic control problem involving a continuum of nonexchangeable and interacting agents and reduces to solving a continuum of Hamilton–Jacobi–Bellman (HJB) equations connected through a graphon that models the connections among agents. Using a fixed-point argument, we prove that this HJB system admits a unique solution in the space of bounded functions when the discount rate is high (i.e., agents are myopic). Under certain assumptions, we also establish regularity properties of the system, such as equi-continuity. We propose a finite difference scheme for computing the HJB system and prove the uniqueness and existence of its numerical solutions. The mean-field logit dynamic is applied to a case study on inland fisheries resource management in the upper Tedori River of Japan. A series of computational cases are then conducted to investigate the dependence of the dynamic on both the discount rate and graphon.

## 1. Introduction

### 1.1 Background

Decision-making among interacting agents drives many social, biological, and ecological phenomena [1-5]. These multi-agent phenomena have been studied across a wide variety of research fields. Two modeling approaches have emerged for such systems: microscopic agent-based models and macroscopic mean-field models, as reviewed below.

Agent-based models track the positions of individual agents in physical and/or phase spaces and have been widely applied to collective motion in both human and non-human animals [6-9], cell movement [10,11], and financial modeling [12,13]. A key advantage of agent-based models lies in their flexibility to encompass complex decision-making rules. However, they often lack mathematical tractability, making rigorous challenges. In contrast, mean-field models provide a statistic description of the collective behavior of a large number of agents. These models offer greater theoretical tractability, which is likely to explain their widespread use in mathematics, physics, and related research fields [14-17]. Mean-field models allow for approximating complex agent-based systems under certain assumptions and extracting their essential features. Recent examples include a theoretical explanation of the macroscopic convergence behavior of genetic algorithms [18], an Ising-type model of the human-AI interactions [19], and an analysis of the stability of socio-environmental dynamics under climate change [20].

One area within mean-field models that is still relatively unexplored is that of graphon games, games involving a large number of possibly heterogeneous (i.e., non-exchangeable) and interacting agents [21]. A graphon, often formulated as the dense limit of a sequence of graphs, serves as an effective modeling tool in a range of applied contexts, including power grids [22,23], epidemic spread [24,25], and machine-learning architecture [26,27]. Both static [28,29] and dynamic graphon games [30-32], along with related models, have been studied. Among these models, graphon mean-field games have attracted attention for their flexibility in describing diverse dynamics shaped by agent heterogeneity, as reviewed below.

A mean-field game is a stochastic control problem driven by noise and played by homogeneous (i.e., exchangeable) agents. It is typically formulated as a system of differential equations consisting of a Hamilton–Jacobi–Bellman (HJB) equation, which characterizes each agent's optimal actions, and a Fokker–Planck (FP) equation, which describes the evolution of the distribution of agent states [33]. Solving a mean-field game yields both the agent's optimal actions, defined as that maximize or minimize a given objective, and the resulting macroscopic behavior of the multi-agent system. Extension of mean-field games to non-Gaussian noise has been studied both theoretically and numerically [34-36]. Applications of mean-field games span a variety of domains, including green finance [37], clean energy management [38], and synchronization phenomena [39].

A graphon mean-field game extends classical ones by incorporating agent heterogeneity and explicitly modeling the network of interactions among agents. Its governing system consists of coupled continua of FP and HJB equations, making it significantly more complex and increasing the degrees of freedom than classical mean-field games [40,41]. Linear-quadratic cases, which admit closed-form solutions, have been well studied [42,43]. More complex cases, such as optimal investment problems involving agent competition, have also been explored both theoretically and numerically [44,45]. The structure of the FP equations under different forms of graphons have also received attention [46]. Although there have been substantial progresses in the analysis of graphon mean-field games, their application to real-world problems remains limited. Specifically, problems concerning the sustainability of natural resources and environmental management, topics of growing global significance, have not been adequately discussed within this framework. This gap is noteworthy given that such management problems inherently involve complex social interactions. These research gaps motivate the present study, whose aims and contributions are explained below.

### 1.2 Aim and contribution

The aim of this study is to formulate, analyze, and compute a novel graphon mean-field game, with a specific application to fisheries resource management. Our contributions toward this aim are as follows:

We investigate mean-field games in a continuous action space, as discrete versions can be obtained by suitably quantizing the actions. Our formulation is motivated by evolutionary game models, a major mean-field framework where agents sequentially update their actions in response to the behavior of others [47]. This setting is naturally suited to resource and energy management problems, where agents aim to optimize a utility that balances the costs and benefits of harvesting a resource through sequential decision-making. We focus on the logit dynamic, where agent actions are chosen via a softmax-type regularized maximization model [48]; logit dynamics are representative mathematical models that describe agents' actions in uncertain environments where the softmax part modulates the size of uncertainties. This dynamic allows for an approximate treatment of Nash equilibria in utility-based games [49]. Variants of the logit dynamic have emerged in a range of domains, including similarity-based learning [50], bilingual game [51], quantal response problem [52], games with Bayesian learning [53,54], and transformer-based architectures in machine learning [55]. However, applications of the logit dynamic to a continuum of heterogeneous agents interacting through a graphon remain largely unexplored.

For homogeneous agents, recent studies have shown that, in the context of evolutionary games, there exists a corresponding mean-field game from where the former can be (approximately) obtained as the most myopic limit of the latter [56-58]. However, these studies lack theoretical results on the well-posedness of the mean-field game. More recently, for a replicator dynamic involving non-exchangeable agents, but without incorporating a graphon structure, a fixed-point approach was proposed to establish the unique existence of solutions to the corresponding mean-field game [59].

In this paper, we use formalism for stationary mean-field games as presented in Chapter 7 of Bensoussan et al. [60] due to its theoretical simplicity. We show that the proposed graphon mean-field game, referred to as the graphon mean-field logit dynamic (G-MFLD), can be analyzed within this approach. The G-MFLD is designed to address dynamic decision-making problems among heterogeneous agents in uncertain environments in a consistent way, i.e., within a single mathematical model by extending the logit formalism. The key elements of this model are utility and graphon, the former to determine agents' preference and latter agents' interactions. The proposed model as a mean-field game enables us to connect the optimized agents' actions and their objectives in a dynamic and feedback (closed-loop) manner, which is natural in social modeling strategies but is absent in classical evolutionary game models. The feedback nature allows for mechanistically understanding social interactions among agents and their consequences. Moreover, the G-MFLD is a well-posed mathematical model. Specifically, we show that the associated HJB system admits a unique solution in a suitable function space of bounded (but not necessarily continuous) functions, at least when the discount rate (i.e., the degree of myopicity) is sufficiently large. The novel contributions of this work are the introduction of the G-MFLD itself and the inclusion of a graphon structure into the game. Particularly, we use a key structural feature of the G-MFLD; the FP and HJB equations can be combined into a single HJB system, to which a fixed-point approach applies. Notably, our formulation is new even in the absence of graphon or agent heterogeneity.

We present a finite difference scheme for computing the HJB system and provide theoretical results on its convergence. The convergence of numerical solutions is achieved when the coefficients involved in the problem are sufficiently regular. Additionally, we briefly investigate convergence in the viscosity sense [61], as viscosity solutions provide a natural framework for addressing optimality equations in stochastic control. Our HJB equation fits within this framework under certain restrictive conditions.

As a computational application of the G-MFLD, we investigate an inland fisheries management problem in a mountainous river system in Japan, which includes a protected area inhabited by the spotless charr, a rare fish species. Although inland fisheries are smaller in scale than marine fisheries, they are recognized as important contributors to sustainable development goals by providing significant economic, recreational, and health benefits [62,63]. Our model represents anglers at each location as agents, each parameterized by their preferred fishing positions within the river system. We computationally investigate the importance of establishing protected areas to conserve rare inland fish species. Consequently, this paper contributes to modeling, computation, and application of a novel mathematical model.

As reviewed above, logit dynamics have been employed as mathematical models for homogeneous agents in uncertain environments. By contrast, our model, particularly the model in the application, accounts for heterogeneity in agents corresponding to fishing spot preferences from a dynamic viewpoint by explicitly linking agents' behavior and their preferences by using graphons, leading to their dynamically optimized actions. Graphons have been proposed as a novel approach to efficiently link agents within the framework of logit dynamics, which theoretically allows for a variety of linkage types ranging from discrete to continuous. The G-MFLD is therefore a versatile framework that provides the information about agents' optimal actions in uncertain environments along with their mechanisms. Such information cannot be obtained by using static or homogeneous models. Our numerical computations are conducted under a simplified condition (e.g., symmetry assumption of graphon) and are exploratory, but we believe that the results obtained highlights the significance of modeling the heterogeneity among agents in applied problems.

The rest of this paper is organized as follows. **Section 2** provides a brief review of the classical logit dynamic. **Section 3** introduces the G-MFLD, derives its HJB system, and presents theoretical results on well-posedness, regularity, and numerical discretization of the system. **Section 4** applies the G-MFLD to a resource management problem. **Section 5** concludes the study and discusses future research directions. The **Appendix** contains the proofs of the propositions presented in the main text.

## 2. Logit dynamic

In this section, we review the classical logit dynamic for homogeneous agents. The discussion is kept at a formal level, while the detailed mathematical framework utilized throughout this paper is introduced in **Section 3**.

### 2.1 Formulation

Let $t \geq 0$ denote time. The space of agent actions is defined as the compact domain $\Omega = [0,1]$. The space of probability measures on $\Omega$ is denoted by $\mathcal{P}(\Omega)$. Given a utility $U : \Omega \times \mathcal{P} \to \mathbb{R}$ to be maximized by choosing agent actions, the logit dynamic is formally described by an evolution equation, which governs the time-dependent probability measure $\mu$ over agent actions [49]:

$$\frac{\mathrm{d}\mu_t(\mathrm{d}x)}{\mathrm{d}t} = \mathfrak{L}_\eta[\mu_t](x)\mathrm{d}x - \mu_t(\mathrm{d}x) \quad \text{for all } x \in \Omega \text{ and } t > 0 \tag{1}$$

with $\mu_0 \in \mathcal{P}(\Omega)$. Here, $\mathfrak{L}_\eta[\mu_t](x)$ represents the logit function

$$\mathfrak{L}_\eta[\mu_t](x) = \frac{e^{\eta U(x,\mu_t)}}{\int_\Omega e^{\eta U(z,\mu_t)}\mathrm{d}z} \tag{2}$$

with a regularization parameter $\eta > 0$.

### 2.2 Logit function

The coefficient $\mathfrak{L}_\eta[\mu]$ in (2) is a softmax function with $\eta^{-1}$ serving as a temperature parameter. A larger value of $\eta$ corresponds to higher temperature, under which $\mathfrak{L}_\eta(\mu)$ as a function of $x$ becomes more concentrated around its maximizer. The regularization parameter $\eta$ measures the level of uncertainty in the environment where agents operate. In fact, for each $x \in \Omega$ and $\mu \in \mathcal{P}(\Omega)$, we have (e.g., Eq. (3) in Lahkar et al. [64])

$$\mathfrak{L}_\eta[\mu_t](x) = \underset{\phi(\cdot)\mathrm{d}x \in \mathcal{P}(\Omega)}{\arg\max} \left\{ \underbrace{\int_\Omega U(x,\mu)\phi(x)\mathrm{d}x}_{\text{Objective to be maximized}} - \frac{1}{\eta} \underbrace{\int_\Omega \phi(x)\ln\phi(x)\mathrm{d}x}_{\text{Relative entropy}} \right\}, \tag{3}$$

where $\phi$ is the Radon–Nikodym derivative, and the right-hand side of (3) represents a maximization problem consisting of an objective function (the first term) and a penalty in the form of relative entropy (second term). The maximizing $\phi$ equals $\mathfrak{L}_\eta[\mu]$. Since lower relative entropy implies that the maximizing $\phi$ is closer to the uniform distribution on $\Omega$ ($\phi \equiv 1$), a larger value of $\eta$ allows for higher relative entropy and hence allows $\phi$ to be more concentrated around the maximizer of the utility. Flynn and Sastry [65] called the last term in (3) the "expected entropy costs," reflecting the idea that making precise choices incurs a cost. From this perspective, the logit dynamic in (1) models a sequential action update process where a smaller value of $\eta$ implies larger uncertainties faced by agents.

***Remark 1*** The relationship (3) has been effectively utilized in the formulation of a more sophisticated mathematical model based on a logit dynamic [66,67].

### 2.3 Key properties

Any equilibrium $\mu_\infty \in \mathcal{P}(\Omega)$ of the logit dynamic (1) satisfies

$$\mathfrak{L}_\eta[\mu_\infty](x)\mathrm{d}x = \mu_\infty(\mathrm{d}x) \quad \text{for all} \quad x \in \Omega, \tag{4}$$

showing that every equilibrium admits a density. According to Theorem 3.1 in Lahkar and Riedel [49], such an equilibrium exists if $\mathfrak{L}_\eta[\mu]$ is continuous with respect to $\mu \in \mathcal{P}(\Omega)$ and $x \in \Omega$. Furthermore, by Theorem 3.4 in Lahkar and Riedel [49], the logit dynamic (1) admits a unique solution $\mu_t \in \mathcal{P}(\Omega)$ ($t \geq 0$) if $U$ is sufficiently regular. However, even under such conditions, the equilibrium may depend on the initial condition $\mu_0$ as demonstrated in Example 5.1 in Lahkar et al. [68]. Another remarkable property of the logit dynamic (1) is that its equilibrium (4) approximates a Nash equilibrium[1] as $\eta$ decreases, under certain conditions on $U$ (Theorem 3.2 [49]).

A key property of equation (1) is that it is a nonlinear FP equation. Specifically, we have

$$\begin{aligned}\underbrace{\frac{\mathrm{d}\mu_t(\mathrm{d}x)}{\mathrm{d}t}}_{\text{Time evolution}} &= \underbrace{\left(\int_{z\in\Omega}\mu_t(\mathrm{d}z)\right)}_{\text{Inflow rate}}\cdot\underbrace{\mathfrak{L}_\eta[\mu_t](x)\mathrm{d}x}_{\text{Inflow probability}} \\ &\quad -\underbrace{\left(\int_{z\in\Omega}\mathfrak{L}_\eta[\mu_t](z)\mathrm{d}z\right)}_{\text{Outflow rate}}\cdot\underbrace{\mu_t(\mathrm{d}x)}_{\text{Outflow probability}}\end{aligned} \quad \text{for all} \quad x \in \Omega \quad \text{and} \quad t > 0. \tag{5}$$

This re-interpretation is the key to the modeling of the G-MFLD, as discussed in **Section 3**.

## 3. Graphon mean-field logit dynamic

We now formulate the G-MFLD along with mathematical and numerical analyses.

### 3.1 Preliminaries and utility

#### 3.1.1 Notations

[1] For the definition of Nash equilibrium ( A Nash equilibrium is a probability measure such that no agent could benefit from changing its strategy if those of all others are fixed. See Section 2 [64].

The space of agent types, representing heterogeneity, is denoted by $I = [0,1]$, while the space of agent actions is again denoted by $\Omega = [0,1]$. We set the joint domain as $\Xi = \Omega \times I$. The space of bounded functions on $\Xi$ is denoted by $B(\Xi)$, which is equipped with the supremum norm $\|f\|_\infty = \sup_{(x,y)\in\Xi} |f(x,y)|$ for any $f \in B(\Xi)$. With this norm, $B(\Xi)$ forms a Banach space (Chapter 3.3 [69]). The space of finite signed measures on $\Omega$ is denoted by $\mathcal{M}(\Omega)$, equipped with the total variation norm $\|\mu\|_{\mathrm{TV}} = \sup_{|f|\le 1} \left| \int_\Omega f(x)\mu(\mathrm{d}x) \right|$ for any $\mu \in \mathcal{M}(\Omega)$, where $f : \Omega \to \mathbb{R}$ is measurable. A continuum of measures on $\Omega$, parameterized by $I$, is denoted by $\{\mu(y,\mathrm{d}x)\}_{y\in I}$. The space of such parameterized measures is denoted by $\mathcal{M}^{(I)}$, which is equipped with the norm $\|\mu\|_{\mathrm{TV}}^{(I)} = \sup_{y\in I} \|\mu(y,\cdot)\|_{\mathrm{TV}}$. The space of probability measures on $\Omega$ parameterized by $I$ is denoted as $\mathcal{P}^{(I)} \left(\subset \mathcal{M}^{(I)}\right)$.

### 3.1.2 Utility

The utility $U : \Omega \times \mathcal{I} \times \mathcal{P}^{(I)} \to \mathbb{R}$ satisfies **Assumption 1**, which will be used to prove the unique existence of solutions to the HJB system.

***Assumption 1***

***(Boundedness)***
$$0 \le U(x,y,\mu) \le \bar{U} \tag{6}$$
*with a constant* $\bar{U} > 0$ *for all* $(x,y,\mu) \in \Omega \times I \times \mathcal{P}^{(I)}$*, and*

***(Lipschitz continuity)***
$$|U(x,y,\mu) - U(x,y,\nu)| \le L_U \|\mu - \nu\|_{\mathrm{TV}}^{(I)} \tag{7}$$
*with a constant* $L_U > 0$ *for all* $(x,y,\mu,\nu) \in \Omega \times I \times \left(\mathcal{P}^{(I)}\right)^2$.

We also use the following **Assumptions 2 and 3**.

***Assumption 2***

*In addition to* ***Assumption 1****, assume the existence of a non-negative, increasing, and continuous function* $\phi_U : \Omega \to [0,+\infty)$ *such that*

***(Equi-continuity)***
$$|U(x_1,y,\mu) - U(x_2,y,\mu)| \le \phi_U(|x_1 - x_2|) \tag{8}$$
*for all* $(x_1,x_2,y,\mu) \in \Omega^2 \times I \times \mathcal{P}^{(I)}$.

***Assumption 3***

*In addition to* ***Assumption 1****, assume the existence of a non-negative, increasing, and continuous function* $\psi_U : I \to [0,+\infty)$ *such that*

***(Equi-continuity)*** $$\left|U(x,y_1,\mu)-U(x,y_2,\mu)\right| \le \psi_U\left(\left|y_1-y_2\right|\right) \quad (9)$$

*for all* $(x,y_1,y_2,\mu) \in \Omega \times I^2 \times \mathcal{P}^{(I)}$.

Homogeneous version of **Assumption 1** is standard in dynamic games and is fundamental to the well-posedness of their associated evolution equations [e.g., 49, 64], and the simplest example would be $U(x,y,\mu)=\int_{\Omega} f(x-z,y)\mathrm{d}z$ with a suitable bivariate function $f$. **Assumptions 2 and 3** are employed for technical reasons to well-pose our model but are considered not so restrictive unless agents do not have too extreme preferences.

### 3.1.3 Graphon

Originally, a graphon is defined as the continuum limit of adjacency matrices from a sequence of graphs [21]. In this paper, however, we do not employ this bottom-up approach; instead, we consider the continuum limit from the outset. Throughout this paper, we assume that each graphon $W : I^2 \to [0,+\infty)$ satisfies the following **Assumption 4**.

***Assumption 4***

***(Symmetry)*** $$W(y,w)=W(w,y) \text{ for all } y,w \in I \quad (10)$$

*and*

***(Integrability)*** $$0 \le \int_{v\in I} W(y,w)\mathrm{d}w \le \bar{W} \text{ for all } y \in I \quad (11)$$

*with a constant* $\bar{W} > 0$.

Symmetry means that the influence between any two agents is mutual. Integrability is a technical assumption that is satisfied, for example, when the graphon is bounded on $I^2$. **Assumption 4** means that agents' interactions are symmetric and not extremely large on average, which is satisfied for constant, piecewise-constant, and uniformly continuous graphons.

### 3.2 Agent dynamics

We explore a mean-field game where agents make decisions based on logit-like dynamics. To this end, and motivated by the logit dynamic (1), we consider the following stochastic differential equation parameterized by $y \in I$, to describe the dynamics of a representative agent [57,58]:

$$\mathrm{d}X_t^{(y)} = \mathrm{d}J_t^{(y)} \text{ for all } t > 0 \quad (12)$$

with initial condition $X_0^{(y)} \in \Omega$. Here, $\left(J_t^{(y)}\right)_{t\geq 0}$ is a càdlàg point process parameterized by $y \in I$, with jump size $z - X_{t-}^{(y)}$, meaning the post-jump position is $z$, where $z$ is generated by a jump measure which is here a possibly time-dependent probability measure $\nu^{(y)} \in \mathcal{P}^{(I)}$. Under this setting, we recover the logit dynamic in the form of the FP equation (5) by setting $\nu^{(y)} = \mathfrak{L}_\eta[\mu_{t-}]\left(X_{t-}^{(y)}\right)\mathrm{d}x$ and omitting the $y$ dependence. We formulate a mean-field game where $\nu^{(y)}$ serves as a control variable, ensuring consistency with a logit dynamic that incorporates agent heterogeneity under certain conditions.

### 3.3 Control formulation

We consider an infinite-horizon maximization problem, where the control is $\nu^{(y)}$:

$$\phi\left(x, y; \nu^{(y)}\right) = \mathbb{E}\left[\int_0^{+\infty} e^{-\delta_y s}\left\{\underbrace{\delta_y U_s\left(X_s^{(y)}, y, m\right)}_{\text{Scaled utility}} - \underbrace{\frac{1}{\eta_y}\int_\Omega p_s^{(y)}(z)\ln p_s^{(y)}(z)\mathrm{d}z}_{\text{Cost}}\right\}\mathrm{d}s \,\middle|\, X_0^{(y)} = x\right] \tag{13}$$

with discount rate $\delta_y > 0$ and entropic parameter $\eta_y > 0$, both depending on $y \in I$ (subscripts are omitted when there is no risk of confusion) for all $(x, y) \in \Omega \times I$, where the time-dependent probability measure $\left(\nu_t^{(y)}\right)_{t\geq 0}$ is assumed to admit a density $\left(p_t^{(y)}\right)_{t\geq 0}$. For the probability measure $\left(\mu_t^{(y)}\right)_{t\geq 0}$ of agent actions, we consider the weighted time average

$$m(y, \mathrm{d}x) = \delta_y \int_0^{+\infty} e^{-\delta_y s}\mu_s(y, \mathrm{d}x)\mathrm{d}s\,. \tag{14}$$

The integrand in (13) resembles (3); the cost represents the relative entropy of the probability density of jump locations, and hence agent actions, meaning that a more concentrated distribution would incur a higher cost. This implies that taking more sophisticated (and hence more effective) strategies for updating agent actions is more costly. Moreover, the cost increases as $\eta$ decreases, corresponding to the logit function deviating from a uniform distribution, as observed in the logit equilibrium. Therefore, the specific form of the cost in our problem is phenomenologically related to the logit dynamic.

Finally, the value function $\Phi : \Xi \to \mathbb{R}$ is the optimized objective function:

$$\Phi(x, y) = \sup_{\nu^{(y)}} \phi\left(x, y; \nu^{(y)}\right) \quad \text{for all } (x, y) \in \Xi\,. \tag{15}$$

The goal of the mean-field game is to find the maximizer of (15) along with the corresponding system dynamics under this optimal control.

### 3.4 HJB system

Following the dynamic programming argument (Chapter 7 [60]), we formally obtain that the value function satisfies a HJB equation:

$$\delta_y \Phi(x,y) = \delta_y U(x,y,m) + \sup_{p(\cdot)\geq 0, \int_\Omega p(x)\mathrm{d}x=1} \left( \int_\Omega p(z)\left(\Phi(z,y)-\Phi(x,y)\right)\mathrm{d}z - \frac{1}{\eta_y}\int_\Omega p(z)\ln p(z)\mathrm{d}z \right) \quad (16)$$

for all $(x,y)\in\Xi$, and the probability measure $\mu$ of agent actions satisfies the FP equation

$$\begin{aligned}\frac{\mathrm{d}\mu_t(y,\mathrm{d}x)}{\mathrm{d}t} &= \left(\int_{z\in\Omega}\mu_t(y,\mathrm{d}z)\right)\cdot p^*(y,x)\mathrm{d}x - \left(\int_{z\in\Omega}p^*(y,z)\mathrm{d}z\right)\cdot\mu_t(y,\mathrm{d}z) \\ &= p^*(y,x)\mathrm{d}x - \mu_t(y,\mathrm{d}x)\end{aligned} \quad (17)$$

for all $(x,y,t)\in\Xi\times(0,+\infty)$. Here, $p^*(y,x)$ is the maximizing $p$ in (16):

$$p^*(y,\cdot) = \underset{p(\cdot)\geq 0, \int_\Omega p(x)\mathrm{d}x=1}{\arg\max} \left( \int_\Omega p(z)\left(\Phi(z,y)-\Phi(x,y)\right)\mathrm{d}z - \frac{1}{\eta_y}\int_\Omega p(z)\ln p(z)\mathrm{d}z \right) = \frac{e^{\eta_y\Phi(\cdot,y)}}{\int_\Omega e^{\eta_y\Phi(z,y)}\mathrm{d}z}, \quad (18)$$

with which we obtain

$$\left( \int_\Omega p(z)\left(\Phi(z,y)-\Phi(x,y)\right)\mathrm{d}z - \frac{1}{\eta_y}\int_\Omega p(z)\ln p(z)\mathrm{d}z \right)_{p=p^*} = \frac{1}{\eta_y}\ln\left(\int_\Omega e^{\eta_y\{\Phi(z,y)-\Phi(x,y)\}}\mathrm{d}z\right). \quad (19)$$

By (18)-(19) along with $m$ in (14), the HJB and FP equations reduce to

$$\delta_y \Phi(x,y) = \delta_y U(x,y,m) + \frac{1}{\eta_y}\ln\left(\int_\Omega e^{\eta_y\{\Phi(z,y)-\Phi(x,y)\}}\mathrm{d}z\right) \quad (20)$$

and

$$\delta_y\left(m(y,\mathrm{d}x) - \mu_0(y,\mathrm{d}x)\right) = \frac{e^{\eta_y\Phi(x,y)}}{\int_\Omega e^{\eta_y\Phi(z,y)}\mathrm{d}z}\mathrm{d}x - m(y,\mathrm{d}x), \quad (21)$$

respectively.

The FP equation (21) can be rewritten as:

$$m(y,\mathrm{d}x) = \underbrace{\frac{\delta_y}{\delta_y+1}\mu_0(y,\mathrm{d}x)}_{\text{Initial condition part}} + \underbrace{\frac{1}{\delta_y+1}\frac{e^{\eta_y\Phi(x,y)}}{\int_\Omega e^{\eta_y\Phi(z,y)}\mathrm{d}z}\mathrm{d}x}_{\text{Logit part}}, \quad (22)$$

which is formally solvable for $m$. The right-hand side of (22) is a weighted mean of two probability measures with $\delta_y$-dependent weights. We obtain $m=\mu_0$ under the large $\delta_y$ provided that $\Phi$ remains bounded independently of $\delta_y$. This boundedness will be confirmed under certain conditions in the next subsection. According to the HJB equation (20), we have $\Phi = U$ in the limit (as proven later).

### 3.5 Mathematical analysis

#### 3.5.1 Well-posedness

We mathematically analyze the unique existence of solutions to the system (20)-(21). Using the relation in (22), we can eliminate $m$ from the system and obtain a single nonlinear integral equation, which we refer to as the HJB system:

$$\Phi(x,y)=U\left(x,y,\frac{\delta_y}{\delta_y+1}\mu_0(y,\mathrm{d}x)+\frac{1}{\delta_y+1}\frac{1}{\int_\Omega e^{\eta_y(\Phi(z,y)-\Phi(x,y))}\mathrm{d}z}\mathrm{d}x\right)+\frac{1}{\delta_y\eta_y}\ln\left(\int_\Omega e^{\eta_y\{\Phi(z,y)-\Phi(x,y)\}}\mathrm{d}z\right)$$
$$\equiv \mathbb{H}_1[\Phi](x,y)+\mathbb{H}_2[\Phi](x,y) \tag{23}$$
$$\equiv \mathbb{H}[\Phi](x,y)$$

for all $(x,y)\in\Xi$. We show that this equation has a unique fixed point in $B(\Xi)$ under certain conditions. Because the FP equation (22) can be easily solved once $\Phi$ is found by solving (23), we focus on the single equation (23).

For later use, we set

$$\bar{\delta}=\sup_{y\in I}\delta_y,\ \underline{\delta}=\inf_{y\in I}\delta_y,\ \bar{\eta}=\sup_{y\in I}\eta_y,\ \underline{\eta}=\inf_{y\in I}\eta_y. \tag{24}$$

Without significant loss of generality, $\delta_y$ and $\eta_y$ are assumed to be bounded on $I$. The following **Proposition 1** presents our first main theoretical result. Its proof is a bit technical and basically uses a Banach fixed-point theorem.

***Proposition 1***

*Under **Assumption 1**, there exists a unique solution $\Phi\in B(\Xi)$ to the HJB system (23) such that $0\le\Phi\le\frac{\bar{\delta}}{\underline{\delta}}\bar{U}$ if*

$$\frac{2\bar{\eta}L_U}{\underline{\delta}+1}e^{2\bar{\eta}\frac{\bar{\delta}}{\underline{\delta}}\bar{U}}+\frac{1}{\underline{\delta}}\left(1+e^{\bar{\eta}\frac{\bar{\delta}}{\underline{\delta}}\bar{U}}\right)\in(0,1). \tag{25}$$

### 3.5.2 Regularity

In **Proposition 1**, we obtained the unique existence of solutions to the HJB system within the uniformly bounded functions. We investigate the regularity of this unique solution under **Assumption 2**. The following **Proposition 2** shows that the continuity of $\Phi$ with respect to the $x$ direction is qualitatively inherited from the continuity of the utility $U$. Notably, the continuity estimate for $\Phi$ is independent of $\eta_y$.

***Proposition 2***

*Assume **Assumption 2**, condition (25), and that $\delta_y>1$ for all $y\in I$. Then, the unique solution $\Phi\in B(\Xi)$ to the HJB system (23) satisfies the following continuity estimate: for all $x_1,x_2\in\Omega$ and $y\in I$,*

$$|\Phi(x_1,y)-\Phi(x_2,y)|\le\frac{\delta_y}{\delta_y-1}\phi_U(|x_1-x_2|). \tag{26}$$

We also investigate continuity in $y$ under **Assumption 3**. For the rest of this sub-section, we assume that $\delta_y$ and $\mu_y$ are constants to isolate the effect of continuity in $\Phi$ due to the regularity of $U$.

***Proposition 3***

*Assume **Assumption 3**, condition (25), and further that $\delta_y = \delta$ and $\eta_y = \eta$ are constants. Then, the unique solution $\Phi \in B(\Xi)$ to the HJB system (23) satisfies the following continuity estimate: for all $x \in \Omega$ and $y_1, y_2 \in I$,*

$$\left|\Phi(x, y_1) - \Phi(x, y_2)\right| \le \left(1 - \frac{1}{\delta}\left(1 + e^{\eta \bar{U}}\right)\right)^{-1} \psi_U\left(\left|y_1 - y_2\right|\right). \tag{27}$$

Combining **Propositions 2-3** yields the following unified continuity result.

***Proposition 4***

*Assume **Assumptions 2-3**, conditions (25), that $\delta_y = \delta > 1$ and $\eta_y = \eta$ are constants. Then, the unique solution $\Phi \in B(\Xi)$ to the HJB system (23) satisfies the following continuity estimate: for all $x_1, x_2 \in \Omega$ and $y_1, y_2 \in I$,*

$$\left|\Phi(x_1, y_1) - \Phi(x_2, y_2)\right| \le \frac{\delta}{\delta - 1}\phi_U\left(\left|x_1 - x_2\right|\right) + \left(1 - \frac{1}{\delta}\left(1 + e^{\eta \bar{U}}\right)\right)^{-1} \psi_U\left(\left|y_1 - y_2\right|\right). \tag{28}$$

### 3.5.3 Large and small discount limits

Unless otherwise specified, we assume a constant $\delta_y = \delta$ case throughout this subsection. We analyze the limits $\delta \to +\infty$ and $\delta \to +0$ of our HJB system, beginning with the limit $\delta \to +\infty$.

***Proposition 5***

*Under the assumptions of **Proposition 4**, as $\delta \to +\infty$, for the solution $\Phi$ to the HJB system (23) it follows that*

$$\Phi(x, y) \to U(x, y, \mu_0) \text{ at all } (x, y) \in \Xi. \tag{29}$$

According to **Proposition 5**, as $\delta \to +\infty$, $\Phi$ approaches the utility $U$, and hence the mean-field game to a logit dynamic.

***Remark 3*** An equi-continuity result analogous to (28) was not available for the replicator model [59] due to the different coupling structures between the FP and HJB equations. The logit model studied in this paper is special because the FP equation can be solved formally, allowing us to obtain the estimate (28).

The other limit $\delta \to +0$ is less trivial. **Propositions 1-4** rely on the largeness of $\delta$; while this assumption is sufficient, it may not be necessary. Therefore, the limit $\delta \to +0$ may still be meaningful depending on the problem setting. The following FP equation provides insights into this issue (the dynamics below are called the discounted logit dynamics):

$$m(y,\mathrm{d}x) = \frac{\delta_y}{\delta_y + 1}\mu_0(y,\mathrm{d}x) + \frac{1}{\delta_y + 1}\frac{e^{\eta_y U(x,y,m)}}{\int_\Omega e^{\eta_y U(z,y,m)}\mathrm{d}z}\mathrm{d}x \tag{30}$$

for all $x, y \in \Xi$, representing a version of (22) where $\Phi$ is replaced by $U$. For constant $\delta_y, \eta_y$, formally taking $\delta \to +0$ in (30) recovers the logit equilibrium (4). This is understandable, as the limit $\delta \to +0$ corresponds to agents considering information over the entire horizon $(0,+\infty)$ and eventually forgetting the initial condition $\mu_0$. By contrast, this would not extend to the G-MFLD. Indeed, we have

$$\delta\{\Phi(x,y) - U(x,y,m)\} = \frac{1}{\eta_y}\ln\left(\int_\Omega e^{\eta_y\{\Phi(z,y)-\Phi(x,y)\}}\mathrm{d}z\right) \tag{31}$$

whose right-hand side formally vanishes as $\delta \to +0$ if we believe the uniform bound $0 \le \Phi \le \bar{U}$ for small $\delta$. This implies that, under $\delta \to +0$, we have $\Phi \equiv \Phi_0 \in [0,\bar{U}]$ with a constant $\Phi_0$. In fact, substituting $\Phi \equiv \Phi_0$ into the right-hand side of (31) results in 0. This formal result suggests that the G-MFLD should be considered as a different dynamic game model from classical logit dynamics.

#### 3.5.4 With-graphon case

We specify the utility $U$, which plays a key role in connecting heterogeneous agents. More specifically, given a function $u : \Omega \times I \times \mathcal{P}^{(I)} \to \mathbb{R}$ satisfying **Assumption 1**, referred to as the local utility, we define the graphon-based utility

$$U(x,y,\mu) = \int_{w \in I} u(x,w,\mu) W(w,y)\mathrm{d}w \tag{32}$$

with a graphon $W$. The relationship (32) implies that the utility is given by a convolution of the local utility with respect to the graphon. Even this simple convolution leads to a complex form of the HJB system. The utility form in (32) was proposed by Caines et al. [41], where it was called a "cost" in the context of a minimization problem. The utility form (32) represents the simplest framework to capture interactions among agents in evolutionary and mean-field games.

If $u$ satisfies **Assumption 1**, then $U$ is bounded; for all $(x,y,\mu) \in \Omega \times I \times \mathcal{P}^{(I)}$,

$$0 \le U(x,y,\mu) \le \int_{w\in I} \bar{U} W(w,y)\mathrm{d}w = \bar{U}\bar{W}. \tag{33}$$

For the Lipschitz continuity, for all $(x,y,\mu,\nu) \in \Omega \times I \times \left(\mathcal{P}^{(I)}\right)^2$, we have

$$|U(x,y,\mu) - U(x,y,\nu)| \le \int_{w\in I} |u(x,w,\mu) - u(x,w,\nu)| W(w,y)\mathrm{d}w \le \bar{W} L_U \|\mu - \nu\|_{\mathrm{TV}}^{(I)}. \tag{34}$$

Moreover, if $u$ satisfies **Assumption 2**, for all $(x_1, x_2, y, \mu) \in \Omega^2 \times I \times \mathcal{P}^{(I)}$, we have

$$\left|U(x_1, y, \mu) - U(x_2, y, \mu)\right| \le \int_{w \in I} \left|u(x_1, w, \mu) - u(x_2, w, \mu)\right| W(w, y) \mathrm{d}w \le \bar{W} \phi_U \left(\left|x_1 - x_2\right|\right). \quad (35)$$

If $u$ satisfies **Assumption 3**, for all $(x, y_1, y_2, \mu) \in \Omega \times I^2 \times \mathcal{P}^{(I)}$, we have

$$\left|U(x, y_1, \mu) - U(x, y_2, \nu)\right| \le \int_{w \in I} u(x, w, \mu) \left|W(w, y_1) - W(w, y_2)\right| \mathrm{d}w \le \bar{U} \int_{w \in I} \left|W(w, y_1) - W(w, y_2)\right| \mathrm{d}w. \quad (36)$$

In summary, **Assumptions 1 and 2** are satisfied by (32) up to a multiplicative constant. If additionally

$$\int_{w \in I} \left|W(w, y_1) - W(w, y_2)\right| \mathrm{d}w \le C_W \psi_U \left(\left|y_1 - y_2\right|\right) \quad (37)$$

for all $y_1, y_2 \in I$ with a constant $C_W > 0$, then **Assumptions 3** is satisfied by the graphon-based utility (32) up to a multiplicative constant.

### 3.6 Numerical discretization

We present a finite difference scheme for computing the HJB system (23), since a closed-form solution does not appear to be available. While the discretization scheme itself is simple and follows previous studies [59], analyzing its convergence requires careful consideration of $\mathbb{H}$.

The domain $\Xi$ is uniformly discretized by grid points $\mathrm{P}_{i,j} = (\hat{x}_i, \hat{y}_j)$. The first and second elements of $\mathrm{P}_{i,j}$ are $\hat{x}_i = (i - 1/2)\Delta x$ ($i = 1, 2, 3, \ldots, N_x$) and $\hat{y}_j = (j - 1/2)\Delta y$ ($j = 1, 2, 3, \ldots, N_y$), respectively. Here, $N_x, N_y \in \mathbb{N}$ are the numbers of vertices in $x$ and $y$ directions, respectively and we set $\Delta x = 1/N_x$ and $\Delta y = 1/N_y$. We assume $N_x \ge 2$. We write the discretized $\Phi$ at $\mathrm{P}_{i,j}$ as $\Phi_{i,j}$. We assume that $\mu_{i,j} = \int_{\Omega_i} \mu_0(y_j, \mathrm{d}x)$ is available beforehand, where $\Omega_i = [(i-1)\Delta x, i\Delta x)$ ($i = 1, 2, 3, \ldots, N_x - 1$) and $\Omega_{N_x} = [1 - \Delta x, 1]$. The discretized $\delta_y, \eta_y$ at $\mathrm{P}_{i,j}$ are denoted as $\delta_j, \eta_j$, respectively.

The graphon-based utility (32) is discretized at $\mathrm{P}_{i,j}$ as

$$U(x, y, m) \to \sum_{l=1}^{N_y} u_{i,l} W_{l,j} \Delta y, \quad (38)$$

where $u_{i,l}$ is the discretization of $u$ (fully determined later to indicate the dependence on $m$) at $\mathrm{P}_{i,l}$ and $W_{l,j} = \frac{1}{\Delta y} \int_{y_l - \Delta y/2}^{y_l + \Delta y/2} W(w, y_j) \mathrm{d}w$. For the HJB system (23), we use the discretization

$$\Phi(x, y) \to \Phi_{i,j} \quad (39)$$

and

$$\frac{1}{\delta_y \eta_y} \ln\left(\int_\Omega e^{\eta_y \{\Phi(z,y) - \Phi(x,y)\}} \mathrm{d}z\right) \to \frac{1}{\delta_j \eta_j} \ln\left(\sum_{k=1}^{N_x} e^{\eta_j \{\Phi_{k,j} - \Phi_{i,j}\}} \Delta x\right). \quad (40)$$

Consequently, we obtain the discretized HJB system

$$\mathbb{G}_{i,j}\left[\Phi^{(n)}\right]=\Phi_{i,j}-\sum_{l=1}^{N_y}u_{i,l}W_{l,j}\Delta y-\frac{1}{\delta_j\eta_j}\ln\left(\sum_{k=1}^{N_x}e^{\eta_j\{\Phi_{k,j}-\Phi_{i,j}\}}\Delta x\right)=0 \quad (41)$$

at all $i=1,2,3,...,N_x$ and $j=1,2,3,...,N_y$.

We solve the system (41) using an iteration method with relaxation (**Algorithm 1**), analogous to the forward Euler one (Section 2.4 [70]). Here, $\Delta t>0$ is the increment of pseudo time, and $\varepsilon>0$ is a small error threshold. The superscript $(n)$ indicates the value at the $n$ th iteration.

***Algorithm 1***

***Step 1.*** *Set an initial guess* $\Phi_{i,j}^{(0)}=0$. *Set* $n=0$.

***Step 2.*** *Compute* $\Phi_{i,j}^{(n+1)}=\Phi_{i,j}^{(n)}-\delta_j\Delta t\mathbb{G}_{i,j}\left[\Phi^{(n)}\right]$ *at all grid points.*

***Step 3.*** *Compute the error* $\varepsilon^{(n)}=\max_{\substack{1\le i\le N_x\\1\le j\le N_y}}\left\{\left|\Phi_{i,j}^{(n+1)}-\Phi_{i,j}^{(n)}\right|\right\}$.

***Step 4.*** *If* $\varepsilon^{(n)}\le\varepsilon$, *then output* $\Phi^{(n+1)}$ *as the numerical solution and terminate the algorithm. If* $\varepsilon^{(n)}>\varepsilon$, *then set* $n\to n+1$ *and go to* ***Step 2***.

The numerical solution obtained from **Algorithm 1** is considered an approximation of the solution to the discretized HJB system when $\varepsilon^{(n)}$ is sufficiently small. The following **Proposition 6** shows that any numerical solution $[\Phi]_{\substack{1\le i\le N_x\\1\le j\le N_y}}$ to the discretized HJB system (41) is uniformly bounded.

***Proposition 6***

*Any numerical solution* $\left[\Phi_{i,j}\right]_{\substack{1\le i\le N_x\\1\le j\le N_y}}$ *to the discretized HJB system (41) satisfies*

$$0\le\Phi_{i,j}\le\bar{U} \text{ at all grid points.} \quad (42)$$

We also establish the unique existence of (numerical) solutions $\left[\Phi_{i,j}\right]_{\substack{1\le i\le N_x\\1\le j\le N_y}}$ to the discretized HJB system (41). The proof is essentially the same with **Proof of Proposition 1**, as both involve bounded functions defined on a compact set and is therefore omitted. The only difference is that, in **Proposition 6**, solutions are considered in the space of functions $\left[\Phi_{i,j}\right]_{\substack{1\le i\le N_x\\1\le j\le N_y}}\in\mathbb{R}^{N_xN_y}$ equipped with the maximum norm $\max_{\substack{1\le i\le N_x\\1\le j\le N_y}}\left|\Phi_{i,j}\right|$.

***Proposition 7***

*Assume that the local utility $u_{i,l}$ as a function of $\left[\Phi_{i,j}\right]_{\substack{1\le i\le N_x\\1\le j\le N_y}}$ satisfies the following discrete Lipschitz continuity analogous to that assumed in **Assumption 1**:*

$$\left|u_{i,j,1}-u_{i,j,2}\right|\le L_U \max_{\substack{1\le i\le N_x\\1\le j\le N_y}} \frac{1}{1+\delta_j}\left|\frac{e^{\eta_j\Phi_{i,j,1}}}{\sum_{k=1}^{N_x} e^{\eta_j\Phi_{k,j,1}}\Delta x}-\frac{e^{\eta_j\Phi_{i,j,2}}}{\sum_{k=1}^{N_x} e^{\eta_j\Phi_{k,j,2}}\Delta x}\right| \tag{43}$$

*for all $1\le i\le N_x$ and $1\le j\le N_y$, where $\left[\Phi_{i,j,1}\right]_{\substack{1\le i\le N_x\\1\le j\le N_y}},\left[\Phi_{i,j,2}\right]_{\substack{1\le i\le N_x\\1\le j\le N_y}}\in\mathbb{R}^{N_xN_y}$, and $u_{i,j,1}$ and $u_{i,j,2}$ represent $u_{i,j}$ with $\Phi_{\cdot,\cdot}=\Phi_{\cdot,\cdot,1}$ and $\Phi_{\cdot,\cdot}=\Phi_{\cdot,\cdot,2}$, respectively. Assume in addition that $\underline{\delta}>0$ is sufficiently large. Then, the discretized HJB system (41) admits a unique solution that satisfies the bound (42).*

To proceed with the convergence analysis of the discretization, we specify the local utility:

$$u(x,w,m)=g\left(x,w,\frac{1}{\delta_w+1}\frac{\int_{q\in\Omega}h(q)e^{\eta_w\Phi(q,w)}\mathrm{d}q}{\int_{z\in\Omega}e^{\eta_w\Phi(z,w)}\mathrm{d}z}\right)\quad\text{at all }(x,w)\in\Xi, \tag{44}$$

where $m$ is given by (22). The functions $g(x,w,\cdot)$ for each $(x,w)\in\Xi$ and $h$ are assumed to be non-negative and Lipschitz continuous, with Lipschitz constants $L_g$ and $L_h$, respectively. We also assume the uniform boundedness condition $0\le g\le\bar{U}$. Terms involving $\mu_0$ are absorbed into $g$ and are assumed to be discretized without any error. We set $\bar{h}=\max_{x\in\Omega}h(x)$. In the simplest case, $g$ would be an affine function.

The discretized version of (44) is set as

$$u(x,y,m)\to u_{i,l}=g\left(x_i,y_l,\frac{1}{\delta_l+1}\frac{\sum_{k=1}^{N_x}h(x_k)e^{\eta_l\Phi_{k,l}}\Delta x}{\sum_{k=1}^{N_x}e^{\eta_l\Phi_{k,l}}\Delta x}\right)\quad\text{at each }\mathrm{P}_{i,j}. \tag{45}$$

Consequently, the discretized HJB system takes the following form: at each $\mathrm{P}_{i,j}$ ($1\le i\le N_x$, $1\le j\le N_y$),

$$\Phi_{i,j}-\sum_{l=1}^{N_y}g\left(x_i,y_l,\frac{1}{\delta_l+1}\frac{\sum_{k=1}^{N_x}h(x_k)e^{\eta_l\Phi_{k,l}}\Delta x}{\sum_{k=1}^{N_x}e^{\eta_l\Phi_{k,l}}\Delta x}\right)W_{l,j}\Delta y-\frac{1}{\delta_j\eta_j}\ln\left(\sum_{k=1}^{N_x}e^{\eta_j\left\{\Phi_{k,j}-\Phi_{i,j}\right\}}\Delta x\right)=0. \tag{46}$$

We now state **Proposition 8**, which concerns the convergence of numerical solutions to the discretized HJB system (41). The condition stated in (47) is satisfied at least if $W$ is continuous on $I^2$.

***Proposition 8***

*Assume that* $\delta_y = \delta$, $\eta_y = \eta$ *are constants. Further, assume* ***Assumptions 2 and 3*** *as well as the assumptions of* ***Propositions 1 and 4****. Also, assume the following condition:*

$$\sup_j \sup_{y \in (y_j - \Delta y/2, y_j + \Delta y/2)} \int_{w \in I} \left| W(w, y_j) - W(w, y) \right| \mathrm{d}w \to +0 , \tag{47}$$

*and that* $N_x = N_y = N$. *Then, for a sufficiently large* $\delta > 0$ *independent of* $N$, *it follows that*

$$\sup_{1 \le i, j \le N} \left| \Phi(x_i, y_j) - \Phi_{i,j} \right| \to +0 \quad \text{as} \quad N \to +\infty . \tag{48}$$

***Remark 3*** Inspecting **Proof of Proposition 8**, especially the inequality (106), suggests that we have the error bound between the true and computed $\Phi$. Indeed, the inequality (106) shows that $\sup_{1 \le i, j \le N} \left| \Phi(x_i, y_j) - \Phi_{i,j} \right| \le C N^{-1}$ with some constant $C > 0$ if the graphon $W$ and functions $g, h$ are sufficiently regular (e.g., they are Lipschitz continuous on their domains). This estimate is in accordance with our computational results (see **Table A1**).

It is natural to ask whether the theory of viscosity solutions [61] applies to our HJB system, given its origin in a stochastic control(-like) problem. Indeed, the HJB system can be interpreted as a degenerate elliptic system associated with a jump-driven process. This viewpoint enables for a deeper investigation into the convergence of numerical solutions toward the unique continuous solution of such a degenerate elliptic equation [71-74]. We present a proposition concerning the viscosity properties of the discretized HJB system. Its result imply that viscosity solution theory may be applicable to the HJB system only in the case of homogeneous agents, i.e., when there is no dependence on $y$.

***Proposition 9***

*For all* $1 \le i \le N_x$ *and* $1 \le j \le N_y$, *the following quantity*

$$\mathbb{S} \equiv -\sum_{l=1}^{N_y} g\left( x_i, y_l, \frac{1}{\delta_l + 1} \frac{\sum_{k=1}^{N_x} h(x_k) e^{\eta_l \Phi_{k,l}} \Delta x}{\sum_{k=1}^{N_x} e^{\eta_l \Phi_{k,l}} \Delta x} \right) W_{l,j} \Delta y - \frac{1}{\delta_j \eta_j} \ln\left( \sum_{k=1}^{N_x} e^{\eta_j \Phi_{k,j}} \Delta x \right) \tag{49}$$

*as a function of* $\left[ \Phi_{i,j} \right]_{\substack{1 \le i \le N_x \\ 1 \le j \le N_y}}$ *is non-increasing for each* $\Phi_{a,j}$ ($1 \le a \le N_x$) *if*

$$\frac{\bar{h}}{\delta_j + 1} L_g \eta_j W_{j,j} \Delta y \le \frac{1}{\delta_j} \quad \text{at all} \quad 1 \le j \le N_y . \tag{50}$$

Moreover, $\mathbb{S}$ *as a function of* $\left[ \Phi_{i,j} \right]_{\substack{1 \le i \le N_x \\ 1 \le j \le N_y}}$ *is not necessarily component-wise non-increasing.*

**Proposition 9** suggests that monotonicity in the sense of Barles and Souganidis [75] fails for the discretized HJB system when the graphon $W$ is non-trivial, i.e., $W_{k,l} \ne 0$ ($k \ne l$). This means that the numerical solutions obtained via the finite difference scheme may not converge in the viscosity sense when graphon-based interactions among agents are present. In contrast, for the homogeneous case $N_y = 1$, convergence in the viscosity sense may still be established using the classical monotonicity based on the approach of Barles and Souganidis [75], e.g., see Theorem 3.9 in Chowdhury et al. [76] and Section 3.2 in Adusumilli [77].

***Remark 4*** In view of **Proposition 9**, the mean-field games with and without agent heterogeneity are essentially different from each other.

***Remark 5*** Imposing a large discount is a technical assumption that is sufficient (but may not be necessary) at least for the infinite-horizon case as studied in this paper because in this case we need to obtain a bounded value function. Moreover, our mathematical analysis relies on the fixed-point argument (e.g., **Proof of Proposition 1**), where an appropriate contraction mapping is obtained with a large discount. From an application standpoint, this assumption will be restrictive when there are many agents with a long-term perspective. Our computation in **Section 4** implies that numerical solutions converge even for small discount that violates the assumptions of the propositions; however, there is no theoretical guarantee that in such a case the G-MFLD admits a unique solution.

## 4. Application

### 4.1 Study target

The target application of this study is a fisheries management problem involving the charr *Salvelinus leucomaenis*, commonly known as *Iwana* in Japan, in the upper Tedori River System, located in the Mount Hakusan, Ishikawa Prefecture, Japan. The study area has been a part of a geopark (*a single, unified geographical areas where sites and landscapes of international geological significance are managed with a holistic concept of protection, education and sustainable development*[2]), called Hakusan Tedorigawa UNESCO Geopark since May 2023. This designation recognizes the region's unique water and geological environment, ecosystem, landscape, and human activities along the Tedori River[3].

Two types of charr inhabit the upstream reaches of the Tedori River: the spotted (typical) charr and spotless (rare) charr (**Figure 1**). The spotted charr is widely distributed across the eastern half of Honshu Island, including the study area [78]. In contrast, the spotless charr is unique to this region and has been designated as a cultural property (natural monument) by Hakusan City, Ishikawa Prefecture since May 30, 1996[4]. In the study area, the spotless (rare) charr inhibit select upper reaches of the river system, while the spotted (common) charr are distributed more broadly throughout the target river system. Fisheries management in this area is authorized by the Hakusan Shiramine Fisheries Cooperative (HSFC), based in Shiramine Village, Hakusan City. According to HSFC, the spotless charr tend to prefer darker habitats and are more difficult to catch than their spotted counterparts and are more difficult to catch[5]. HSFC has been exploring ways to promote revitalization and sustainable tourism in Mount Hakusan to support the livelihoods of residents. The spotless charr, as a rare species requiring conservation, has been playing a key role in these discussions. For example, there is an ongoing discussion about balancing the conservation of the spotless charr with its potential utilization in tourism promotion. Attracting more tourists, such as anglers, could benefit the local economy; however, it is crucial to avoid exceeding the area's tourism carrying capacity [79,80] to ensure the regional ecosystem is preserved and not degraded.

[2] UNESCO Global Geoparks https://www.unesco.org/en/iggp/geoparks/about (Last accessed on June 23, 2025).

[3] About Hakusan Tedorigawa UNESCO Geopark https://hakusan-geo.jp/about/ (In Japanese. Last accessed on June 23, 2025).

[4] Designated cultural property (natural monument) in Hakusan City. https://www.city.hakusan.lg.jp/bunka/bunkazai/1006096/1002339/index.html (In Japanese. Last accessed on June 23, 2025).

Establishing sustainable tourism is a global issue that has been studied using evolutionary game theory [81-83]. However, it has yet to be studied from the standpoint of graphon (mean-field) games. The Mount Hakusan region, part of a biological reserve, has experienced a slower rate of depopulation than similar rural areas such as reserves, likely due to its higher attractiveness to visitors [84]. This suggests that efforts to protect the regional environment and ecosystem may also help mitigate local population decline indirectly. We apply the G-MFLD framework to a simplified version of the fisheries management problem described above, using it as a case study to explore potential strategies for managing charr populations. Although our case study involves simplified conditions, it provides meaningful insights into sustainable fisheries management for the future.

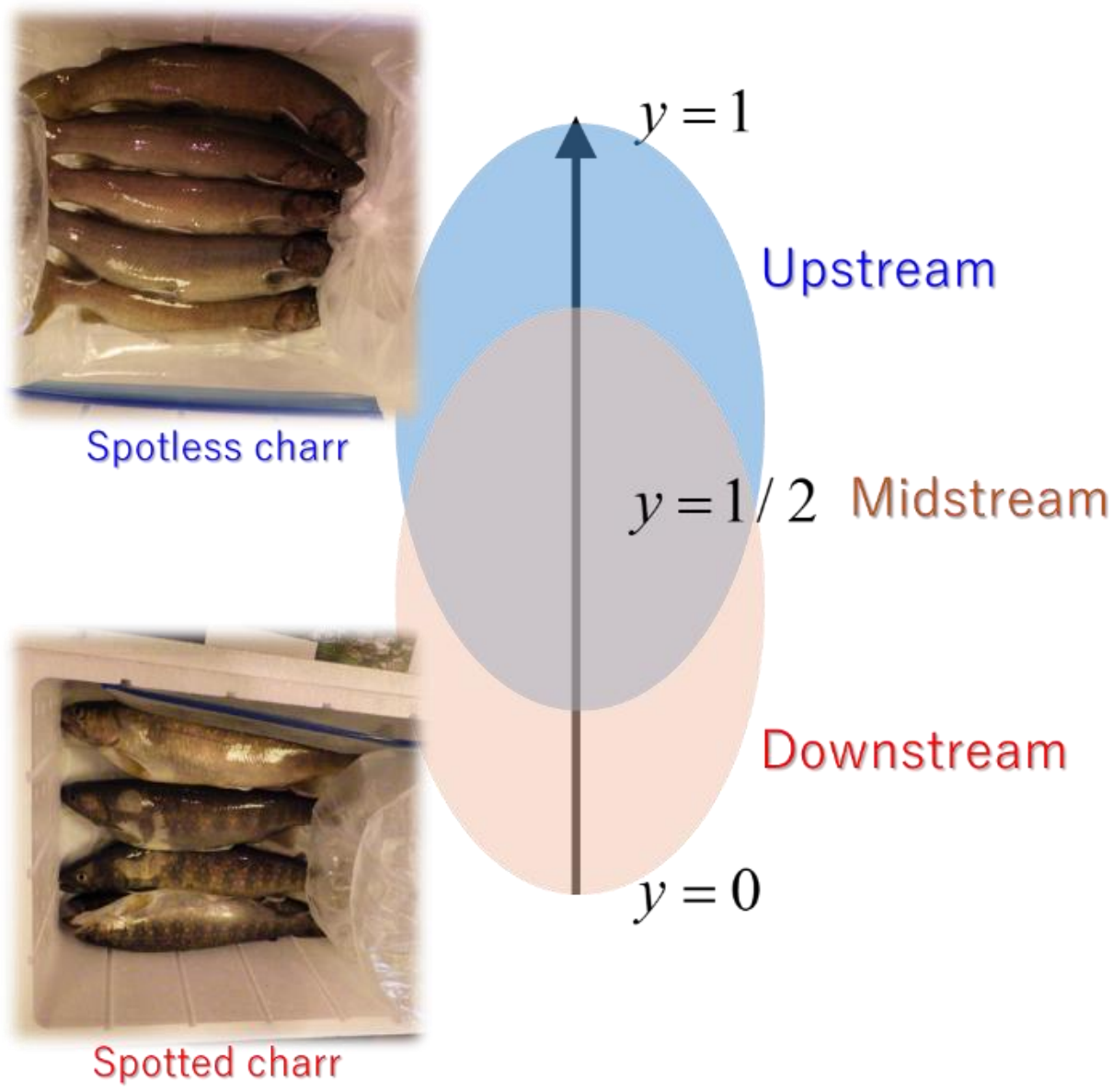

**Figure 1.** Schematic diagram of the application study

[5] Personal communication on May 28, 2025.

### 4.2 Model formulation

Currently, there is no maximum catch limit of the fish in the study site. We therefore investigate how the overfishing of the fish can be suppressed by socially designing the utility. We assume a fisheries management setting along a river reach, where spotless and spotted charr are more abundant in upstream and downstream parts of it, respectively (**Figure 1**). In this context, the agents are anglers. We assume that fishing becomes more costly, requiring greater effort, as one moves farther upstream. Additionally, the catchability of fish at each point in the river decreases as fishing pressure, denoted by $x$ in our context, increases. This captures a density-dependent effect, where higher fishing pressure reduces the overall gain for agents. For simplicity, we model the river reach as the interval $I=[0,1]$ with $y=0,1$ representing the downstream and upstream ends, respectively. We assume spotless charr dominant in $y>1/2$, while spotted charr are dominant in $y<1/2$. The variable $x$ represents fishing pressure, normalized over the interval $\Omega=[0,1]$.

With these assumptions in mind, we specify the graphon-based utility $u$ as follows. This setup corresponds to be one of the simplest forms of potential games with a continuous action space (Section 4 [85]). For all $(x,y)\in\Xi$ and $m\in\mathcal{P}^{(I)}$, we set

$$u(x,y,m)=xA(\alpha_y)-xc(y) \quad \text{with} \quad \alpha_y=\int_{x\in\Omega} xm(y,\mathrm{d}x) \tag{51}$$

with a decreasing function $A:\Omega\to[0,+\infty)$ and a bounded and continuous function. The cost of harvesting at $y$ is $xc(y)$, while the gain is $xA(\alpha_y)$, which decreases as the average fishing pressure $\alpha_y$ increases. The cost $c(y)$ includes both the price of a fishing permit and the travel cost to reach the location. We consider it an increasing function of $y$, since higher fishing fees for rarer species can support conservation efforts, and traveling upstream requires more time. For example, the Tedori River is known for its steep gradient [86], meaning anglers are likely to exert more effort when moving upstream.

As shown in (32), the utility $U$ is defined as a convolution between $u$ and a graphon $W$. In the context of our application, $W$ connects the utilities and hence the decision-making of individual agents across distinct locations $y$. Therefore, each agent's cost and gain reflect both local and global fishing effects. We set $W$ and $c$ as follows:

$$W(y,w)=C_W e^{-\frac{(y-w)^2}{2\theta^2}} \quad \text{for all} \quad y,w\in I \tag{52}$$

and

$$c(y)=c_0+\frac{1}{2}(c_1-c_0)\left(1+\tanh\left(\rho\left(y-\frac{1}{2}\right)\right)\right) \quad \text{for all} \quad y\in I. \tag{53}$$

In (52), $\theta>0$ is a parameter controlling the connection width of the graphon (**Figure 2**) and $C_W>0$ is a normalization constant ensuring that $\int_{y\in I} W(y,\cdot)\mathrm{d}y=1$. In (53), parameters $c_0,c_1>0$ represent the fishing costs for spotless and spotted charr, respectively, and $\rho > 0$ is a parameter controlling the transition of costs across the point $y = 1/2$ (**Figure 3**). This formulation models the charr conservation strategy by dividing the river system into distinct zones: conservation (upstream) and normal fishing regulation (downstream) [78]. The form of graphon (52) implies that the utility at a point in the river is influenced by those around this point, and the parameter $\theta$ controls the width of influence. Under the present setting, a wider graphon width means that each angler exerts fishing pressure in a wider area in the river reach, and his/her target area is more localized as $\theta$ decreases. Intuitively, anglers having a wider width of influence may have more chance of fishing for the rare charr species, possibly resulting in higher fishing pressure in the upstream reach.

We set $A$ as

$$A(\alpha) = 1/\sqrt{\alpha} \quad \text{at all} \quad \alpha > 0, \tag{54}$$

which represents the simplest situation where the gain decreases as the net fishing pressure increases. In the absence of no-graphon settings ($U = u$), $y \in I$ acts as an external parameter. The set of Nash equilibria for the utility $U$ for $y \in I$ is any $m \in \mathcal{P}$ such that $\alpha_y$ maximizes the quasi-potential, defined by (Lemma 1 and Corollary 2 [85])

$$\int_0^{\alpha_y} A(\alpha)\mathrm{d}\alpha - c(y)\alpha_y = 2\sqrt{\alpha_y} - c(y)\alpha_y. \tag{55}$$

Then, an elementary calculation shows that any Nash equilibrium satisfies $\alpha_y = 1/\min\{c(y),1\}^2$. The utility in this case does not satisfy **Assumption 1** but can be made so by regularizing $\sqrt{\alpha}$ as $\sqrt{\alpha+\gamma}$ with $\gamma$ being a small positive constant.

The computational resolution is set to $N_x = N_y = 300$. The time increment for iteration is $\Delta t = 0.01$, with an error threshold of $\varepsilon = 10^{-10}$. Numerical solutions are typically achieved within $O(10^3)$ to $O(10^5)$ iterations for our computational cases. For the cost (53), we use the following parameter values unless otherwise specified, with which the upstream ($y > 1/2$) and downstream reaches ($y < 1/2$) have contrasting optimized $\alpha_y = 1/\{c(y)\}^2$ values and they are separated with a gradual transition: $c = \sqrt{2}$, $c = \sqrt{10}$, and $\rho = 0.05$. We assume a uniform initial condition $\mu_0(y,\mathrm{d}x) = \mathrm{d}x$ ($x \in \Omega$, $y \in I$), representing the most random distribution. In this setting, since $\mu_0$ admits a density, $m$ also does. We therefore investigate the probability density $p$ of $m$ ($p(x,y)\mathrm{d}x = m(y,\mathrm{d}x)$).

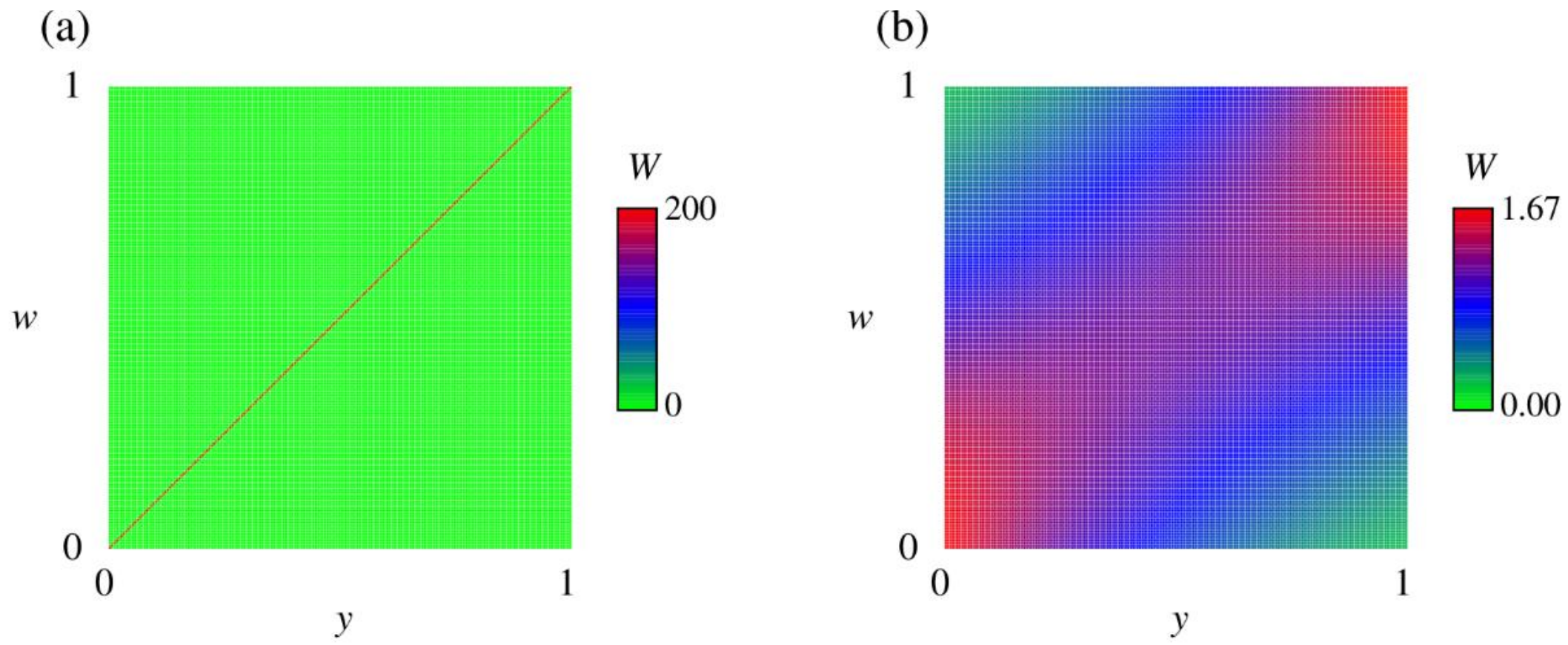

**Figure 2.** Computed graphon $W$ for $\theta$ values of (a) 0.0001 (blue) and (b) 0.5. We set $N_y = 200$

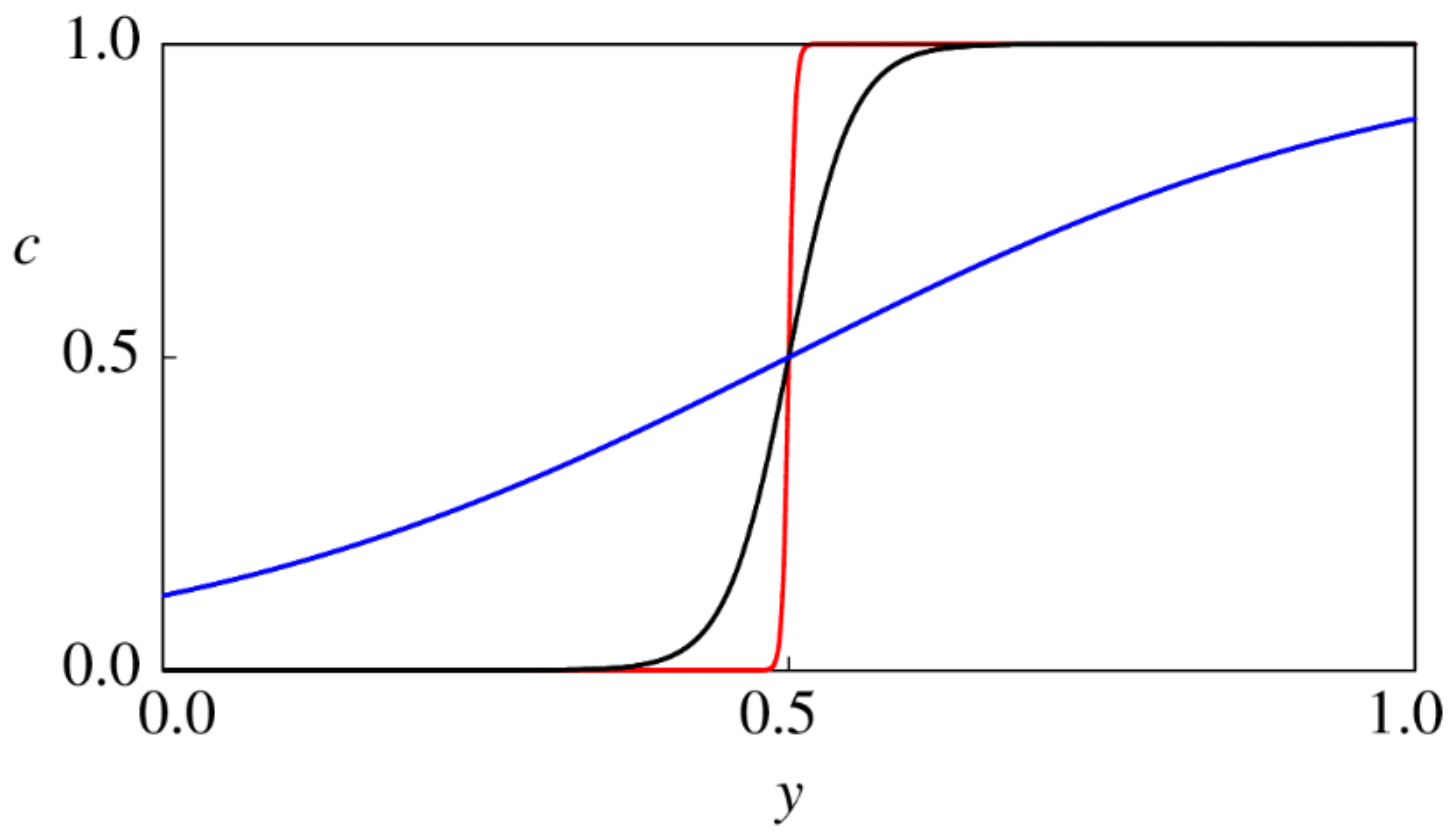

**Figure 3.** Cost $c$ for $\rho$ values of 0.5 (blue), 0.05 (black), and 0.005 (red). We set $c_0 = 0$ and $c_1 = 1$

### 4.3 Results and discussion

#### 4.3.1 No-graphon case

We first study the no-graphon case, corresponding to (52) with $\theta \to +0$ that is computationally addressed by choosing $\theta = 0.0001$. We examine distinct values of the discount rate $\delta$ and the regularization parameter $\eta$, both being constant. We consider four cases covering all combinations of large and small $\delta$ and $\eta$ values: Case (A) $(\delta,\eta) = (0.5, 1/0.5)$, Case (B) $(\delta,\eta) = (0.5, 1/0.005)$, Case (C) $(\delta,\eta) = (0.005, 1/0.5)$, and Case (D) $(\delta,\eta) = (0.005, 1/0.005)$.

**Figures 4 and 5** show the computed distributions of $p$ for Cases (A)-(D) in the discounted logit dynamic (30) and the G-MFLD, respectively. **Figure 6** shows the computed average fishing pressure $\alpha = \alpha_y$ based on $p$ presented in **Figures 4 and 5**. **Figure 4** suggests the role of the parameters $\delta$ and $\eta$ in the discounted logit dynamic; smaller values of $\delta$ leads to more contrasting probability densities between large and small $y$, while larger $\eta$ leads to a sharper transition of the probability densities along the $x$ direction. Similar observations apply to the G-MFLD as shown in **Figure 5**. A difference between the discounted logit dynamic and the G-MFLD is that $p$ in the latter exhibit smoother profiles. This difference is visualized in terms of $\alpha$ as shown in **Figure 6**, showing that the variation of the average as a function of $y$ is sharper for the logit dynamic than the G-MFLD. Another important finding from **Figure 6** is that the average fishing pressure in Case (D) of the discounted logit dynamic closely aligns with the theoretical best-response solution, suggesting that this case provides a reasonably accurate approximation of the limiting behavior $\delta, \eta^{-1} \to +0$ ). For the G-MFLD, specifying a small $\delta$ yields an almost constant $\Phi$, supporting the conjecture made in **Section 3.5.3**.

From a fisheries management perspective, the computational results of the G-MFLD suggest that incorporating future predictions, i.e., using a conditional expectation as the objective function, leads to less variable action profiles. This suggests a certain robustness in angler's preferred actions against perturbations, even though these actions may differ from the Nash equilibria, as shown in **Figure 6**. Therefore, the presented mean-field game model should be considered different from the classical and discounted logit models. The results also suggest that higher cost $c$ lead to lower average fishing pressure $\alpha$, and that upstream locations, with a larger $y$ value results in the smaller fishing pressure. Imposing higher fishing fees in the upstream region ( $y > 1/2$ ) is therefore expected to reduce fishing pressure in those areas. However, the effectiveness of this strategy depends on the degree of myopicity ( $\delta$ ) and the level of environmental uncertainty ( $\eta$ ). Specifically, protecting the upstream reach by incurring higher costs becomes more effective when anglers adopt a longer perspective and when environmental uncertainty is relatively low. Environmental education targeting anglers, especially regarding the conservation of charr, and particularly the spotless ones, could be a valuable approach for promoting sustainable fisheries management. Environmental uncertainty may also be mitigated through monitoring fish stocks within the river reach, conducted by fisheries cooperatives or local government authorities. While such activities are not directly considered in the proposed mathematical model, the implications suggest the importance of developing sustainability-oriented thinking among anglers. In the study area, ecological education aimed at promoting sustainable relationships between humans and nature has recently been initiated by researchers [87]. Extending and continuing such activities, particularly among local communities and anglers, would effectively support sustainable fishing practices in the region.

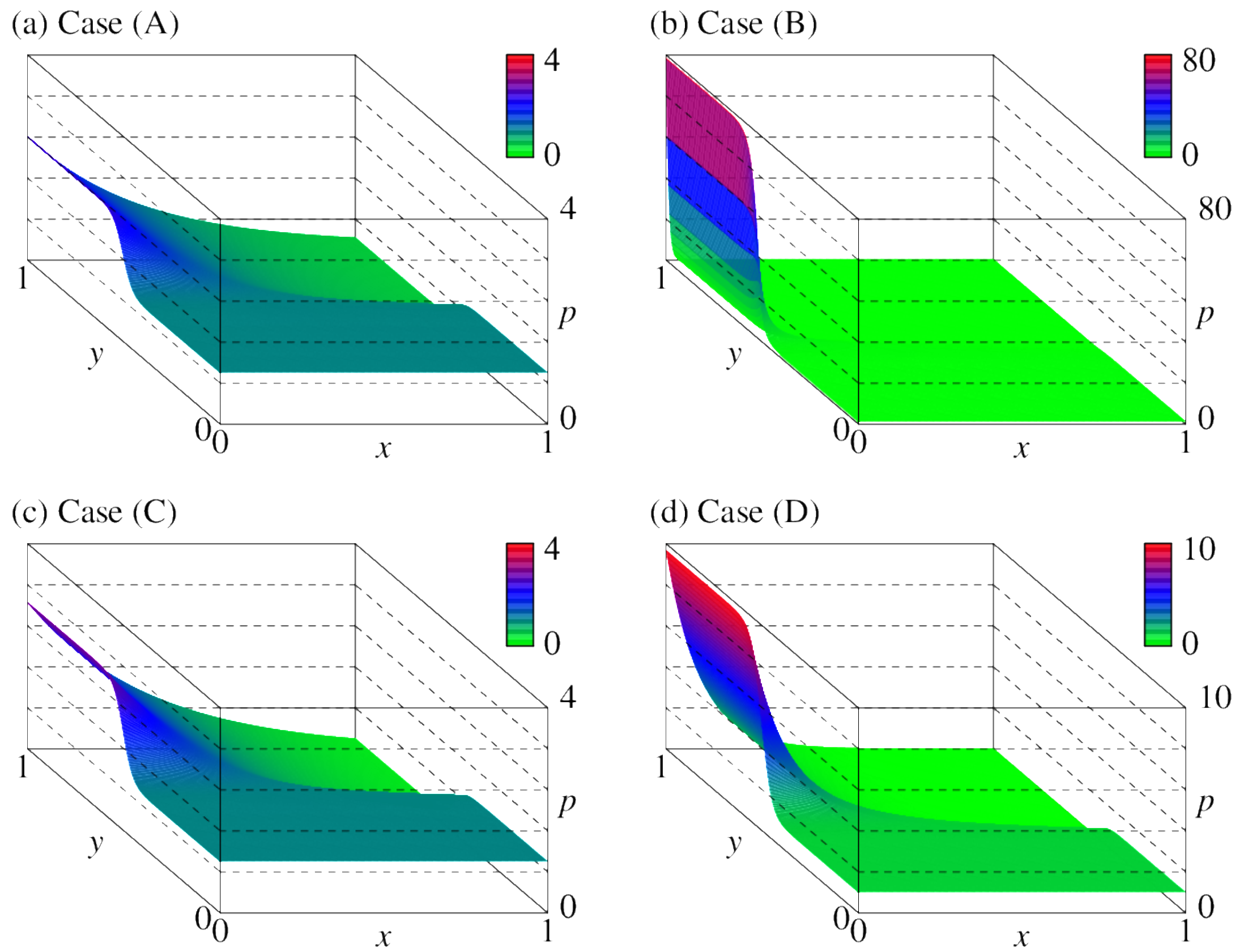

**Figure 4.** Computed probability density $p$ for Cases (A)-(D) for the discounted logit dynamic (No-graphon case): Case (A) $(\delta,\eta)=(0.5,1/0.5)$, Case (B) $(\delta,\eta)=(0.5,1/0.005)$, Case (C) $(\delta,\eta)=(0.005,1/0.5)$, and Case (D) $(\delta,\eta)=(0.005,1/0.005)$.

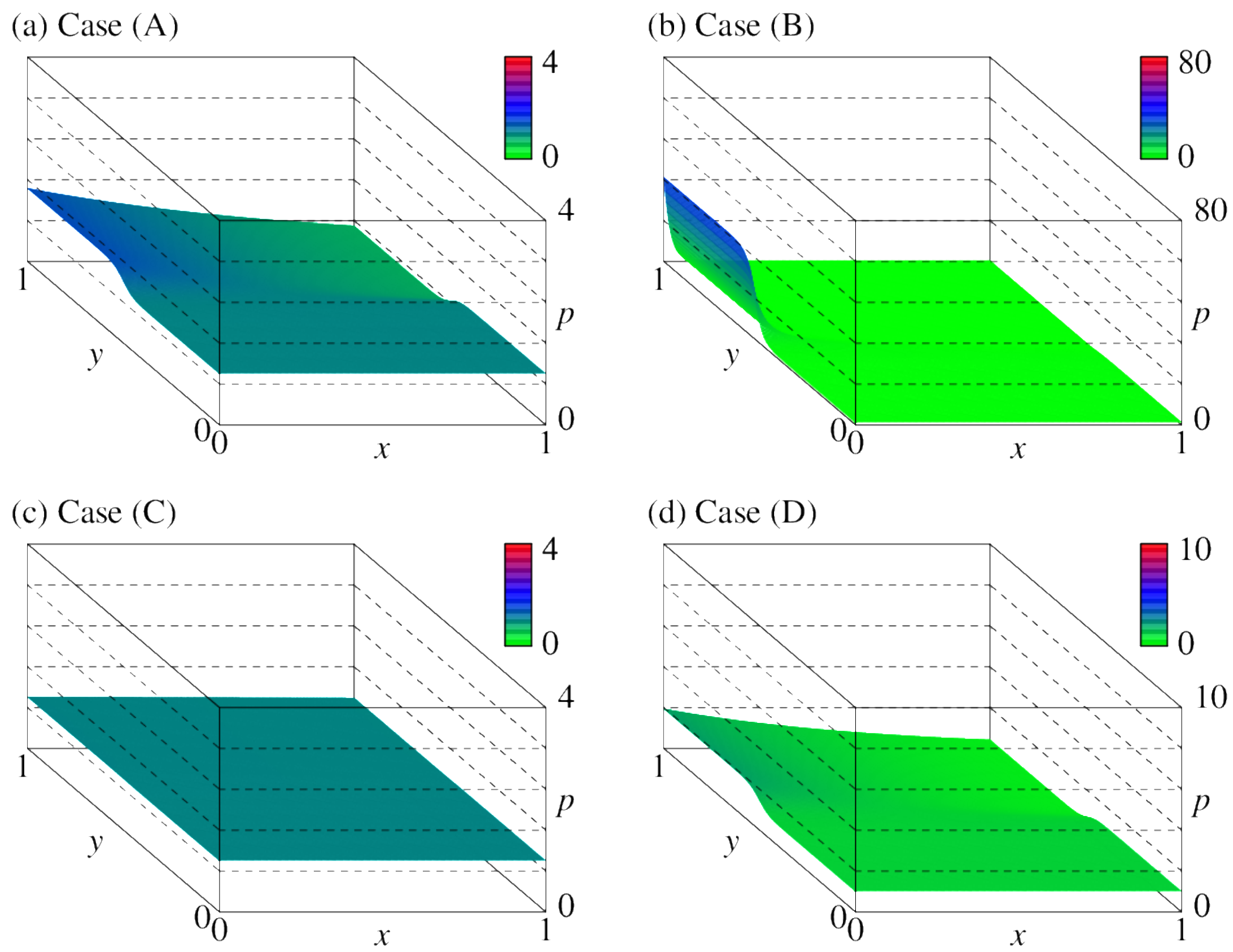

**Figure 5.** Computed probability density $p$ for Cases (A)-(D) for the G-MFLD (No-graphon case): Case (A) $(\delta,\eta)=(0.5,1/0.5)$, Case (B) $(\delta,\eta)=(0.5,1/0.005)$, Case (C) $(\delta,\eta)=(0.005,1/0.5)$, and Case (D) $(\delta,\eta)=(0.005,1/0.005)$.

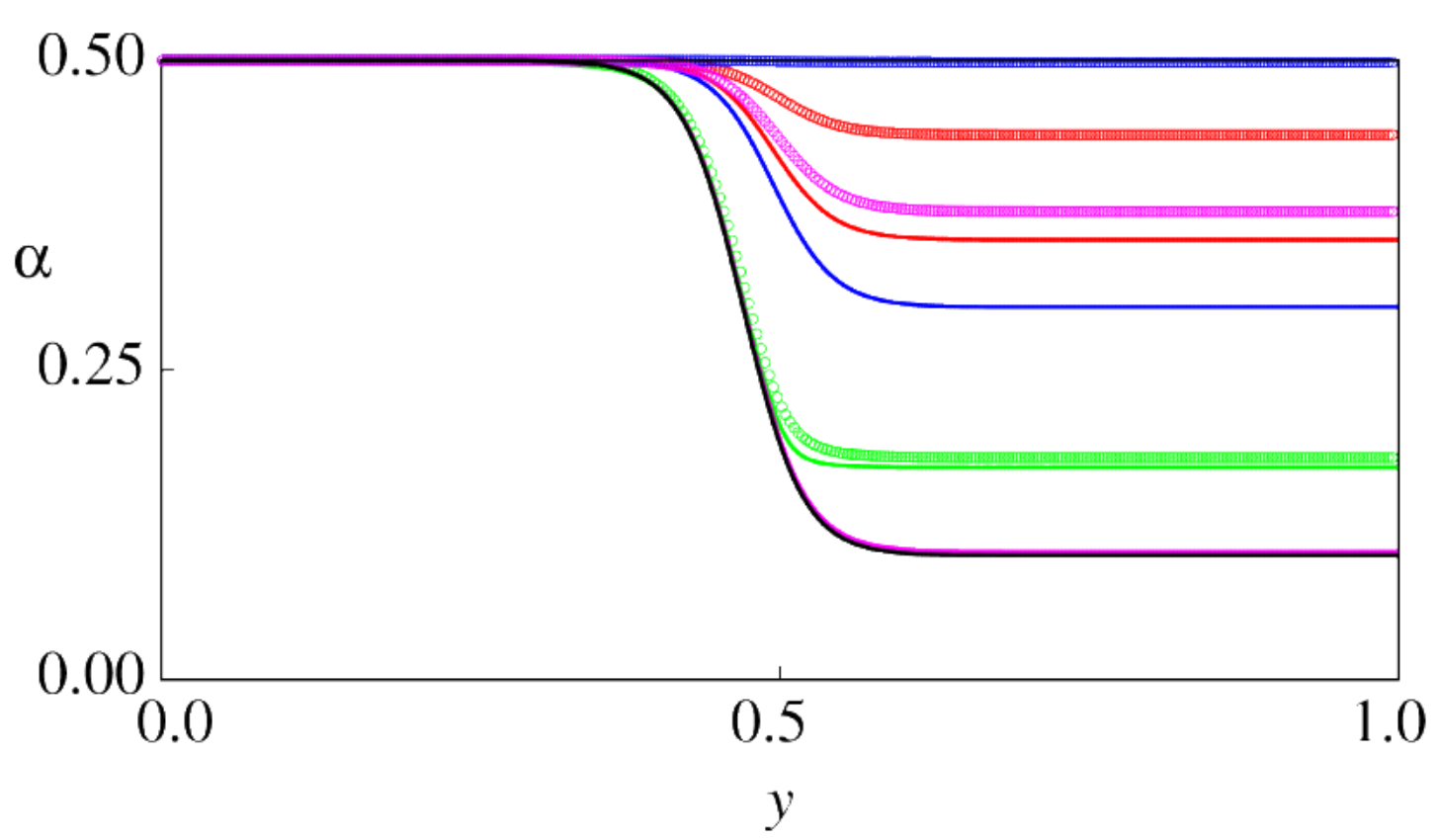

**Figure 6.** Computed average fishing pressure $\alpha=\alpha_y$ based on $p$ presented in **Figures 4 and 5**. Black curve represents the average of Nash equilibria. Colored results correspond to Cases (A)-(D): Case (A) (Red), Case (B) (Green), Case (C) (Blue), and Case (D) (Magenta). The non-black lines correspond to the discounted logit dynamic and non-black circles to the G-MFLD.

#### 4.3.2 With-graphon case

We next consider the with-graphon case, corresponding to (52) with $\theta = 0.5$. **Figures 7 and 8** show the computed distributions of $p$ for Cases (A)-(D) under the discounted logit dynamic and the G-MFLD, respectively. A comparison between **Figures 4 and 7** for the discounted logit dynamic suggests that the presence of the graphon has a regularizing effect in the $y$ direction. This is expected, as the graphon-based utility is a weighted average of local utilities. Similar effects are observed in the G-MFLD as shown by comparing **Figures 5 and 8**.

**Figure 9** shows the computed $\alpha = \alpha_y$ based on $p$ presented in **Figures 7 and 8**. A comparison between the no-graphon case (**Figure 6**) and the with-graphon case (**Figure 9**) implies the significant role of the graphon in shaping $\alpha$, which represents the average fishing pressure at each location in the river reach. While $\alpha$ values decrease in the $y$ direction, similar to the theoretical profile, but exhibit qualitatively different shapes. These profiles are smoother than those in the no-graphon case, and overshoots (i.e., $\alpha$ exceeds $c_0^{-2} = 1/2$) are observed for small $y$ values. These overshoots become more pronounced as $\eta$ decreases.

To better understand the role of graphon, we conduct additional numerical experiments. **Figure 10** shows the computed $\alpha = \alpha_y$ in Case (D) for different values of $\theta$ in the graphon $W$; smaller values of $\theta$ correspond to sharper variations of $W$ and hence weaker connections among agents at distinctive locations. For both the discounted logit dynamic and G-MFLD, smaller values of $\theta$ lead to more gradually varying profiles of $\alpha$. For the discounted logit dynamic, the profile of $\alpha$ is monotone when $\theta = 0.5$ is large but exhibits an overshoot for small $y$. As $\theta$ decreases, the profile of $\alpha$ ($\theta = 2^{-2}, 2^{-3}$) becomes unimodal. These non-monotone profiles are not found in the G-MFLD. This difference likely arises due to the discounted logit dynamic, governed by equation (30), does not originate from a stochastic control framework, unlike the G-MFLD. Moreover, simply integrating the logit dynamic (1) using the relationship in (14) does not yield (30).

Finally, we study cases with spatially variable $\delta_y$, where we assume it as a linear function $0.005 + 0.095(1-y)$ (Case R: longer perspective for rare species) or $0.005 + 0.095y$ (Case M: longer perspective for major species) for all $y \in I$. We fix the constant case $\eta_y = 1/0.005$. **Figure 11** shows the computed average fishing pressure $\alpha = \alpha_y$ for Cases R and M with $\theta = 2^{-1}$ and $\theta = 2^{-7}$. The results in this figure imply that $\alpha_y$ effectively decreases for $y$ in Case R if agents follow the discounted logit dynamic. By contrast, for the G-MFLD, Case R with the narrower graphon width also results in a decreasing $\alpha_y$; however, Case R with the wider graphon width gives a non-monotone $\alpha_y$. A non-monotone $\alpha_y$ is also obtained in Case M with the wider graphon width. Cases M and R with higher costs $(c_0, c_1) = (\sqrt{4}, \sqrt{14})$ are also examined (**Figure 12**), showing that incurring higher cost $c(y)$ lowers $\alpha_y$ without significantly affecting their qualitative profiles. Hence, protecting the rare species would be an effective strategy to protect the rare species for both Cases R and M where the discount rate is distributed in space.

The computational findings in this study imply that the discounted logit dynamic and the G-MFLD represent different mean-field models, particularly when $\theta$ is large. The G-MFLD results show how increasing the effective width of the graphon i.e., $\theta$ influences the spatial distribution of fishing pressure. Specifically, smaller average fishing pressure $\alpha$ is achieved when $\theta$ is small, showing that upstream and downstream reaches should be managed differently. This supports the idea that the upstream region can more effectively function as a protected area for charr. Overall, the results obtained suggest that agents tend to adopt more sustainable fishing behaviors under higher fishing costs. However, when the utility of anglers is not localized, this can result in higher overall fishing pressure. Therefore, areas designated for protection should be subject to significantly higher fishing costs than the other areas to ensure the long-term conservation of rare species.

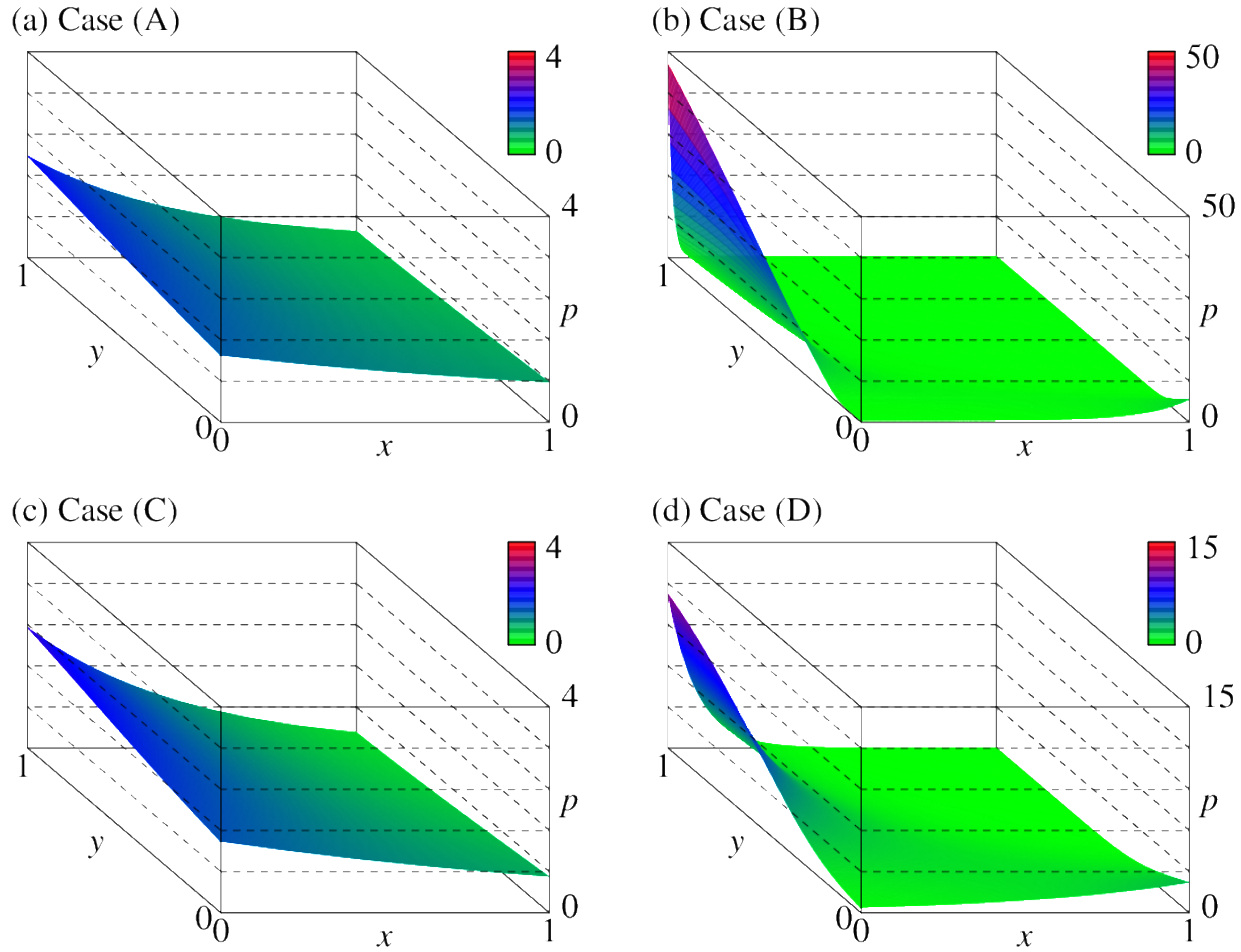

**Figure 7.** Computed probability density $p$ for Cases (A)-(D) for the discounted logit dynamic (With-graphon case): Case (A) $(\delta,\eta)=(0.5,1/0.5)$, Case (B) $(\delta,\eta)=(0.5,1/0.005)$, Case (C) $(\delta,\eta)=(0.005,1/0.5)$, and Case (D) $(\delta,\eta)=(0.005,1/0.005)$.

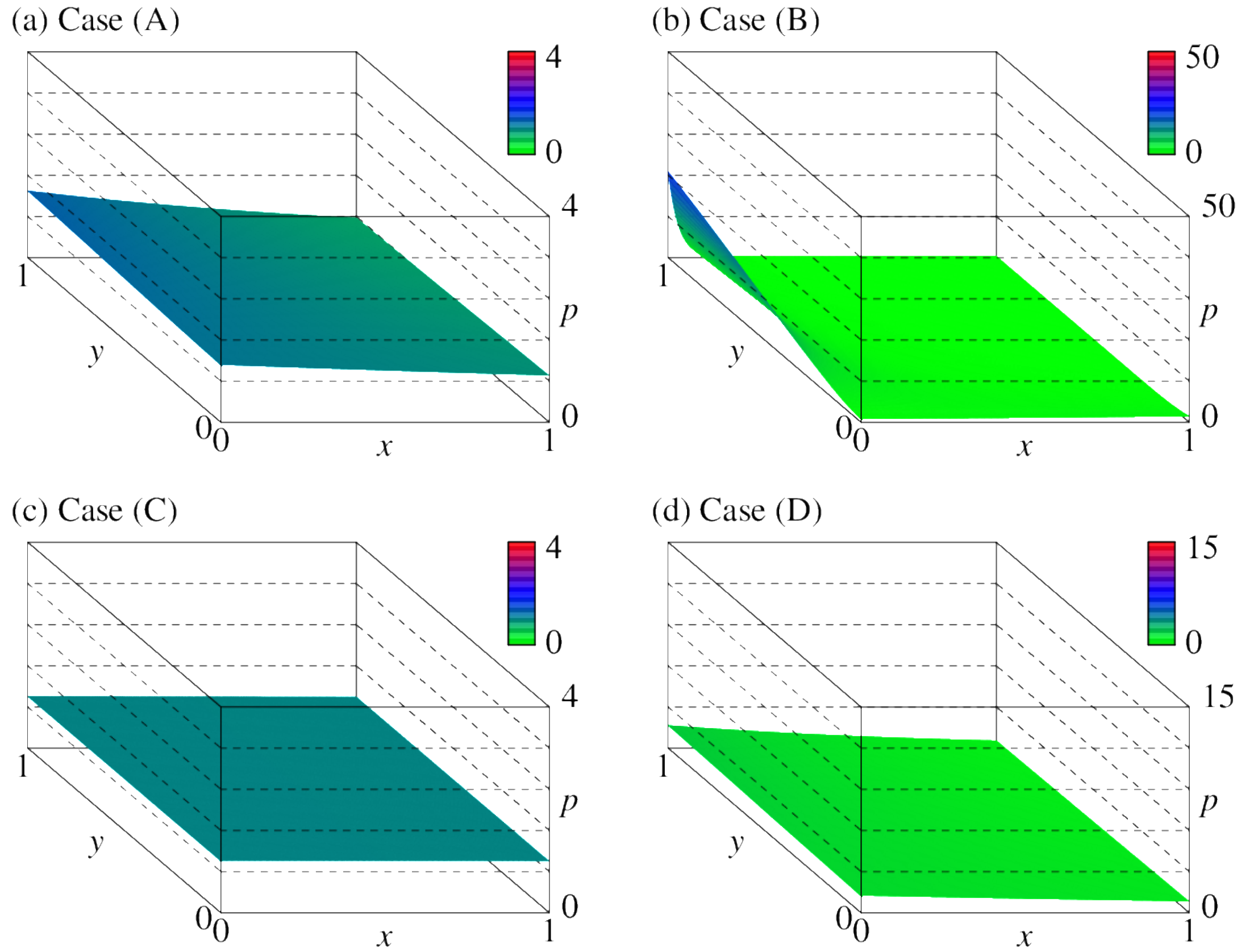

**Figure 8.** Computed probability density $p$ for Cases (A)-(D) for the G-MFLD (With-graphon case): Case (A) $(\delta,\eta)=(0.5,1/0.5)$, Case (B) $(\delta,\eta)=(0.5,1/0.005)$, Case (C) $(\delta,\eta)=(0.005,1/0.5)$, and Case (D) $(\delta,\eta)=(0.005,1/0.005)$.

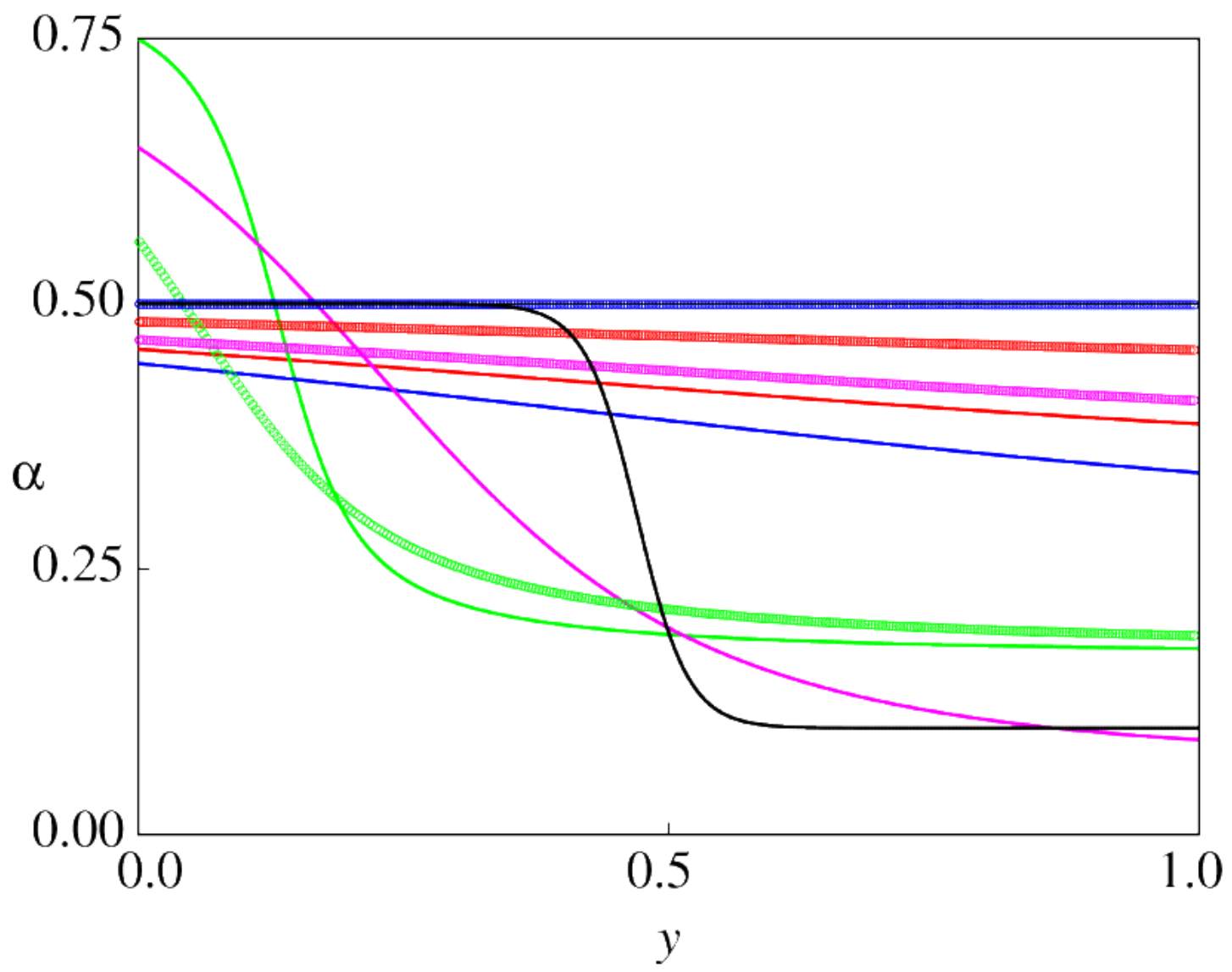

**Figure 9.** Computed average fishing pressure $\alpha=\alpha_y$ based on the probability densities $p$ presented in **Figures 7 and 8**. Black curve represents the average of Nash equilibria. Colored results correspond to Cases (A)-(D): Case (A) (Red), Case (B) (Green), Case (C) (Blue), and Case (D) (Magenta). The non-black lines correspond to the discounted logit dynamic and non-black circles to the G-MFLD.

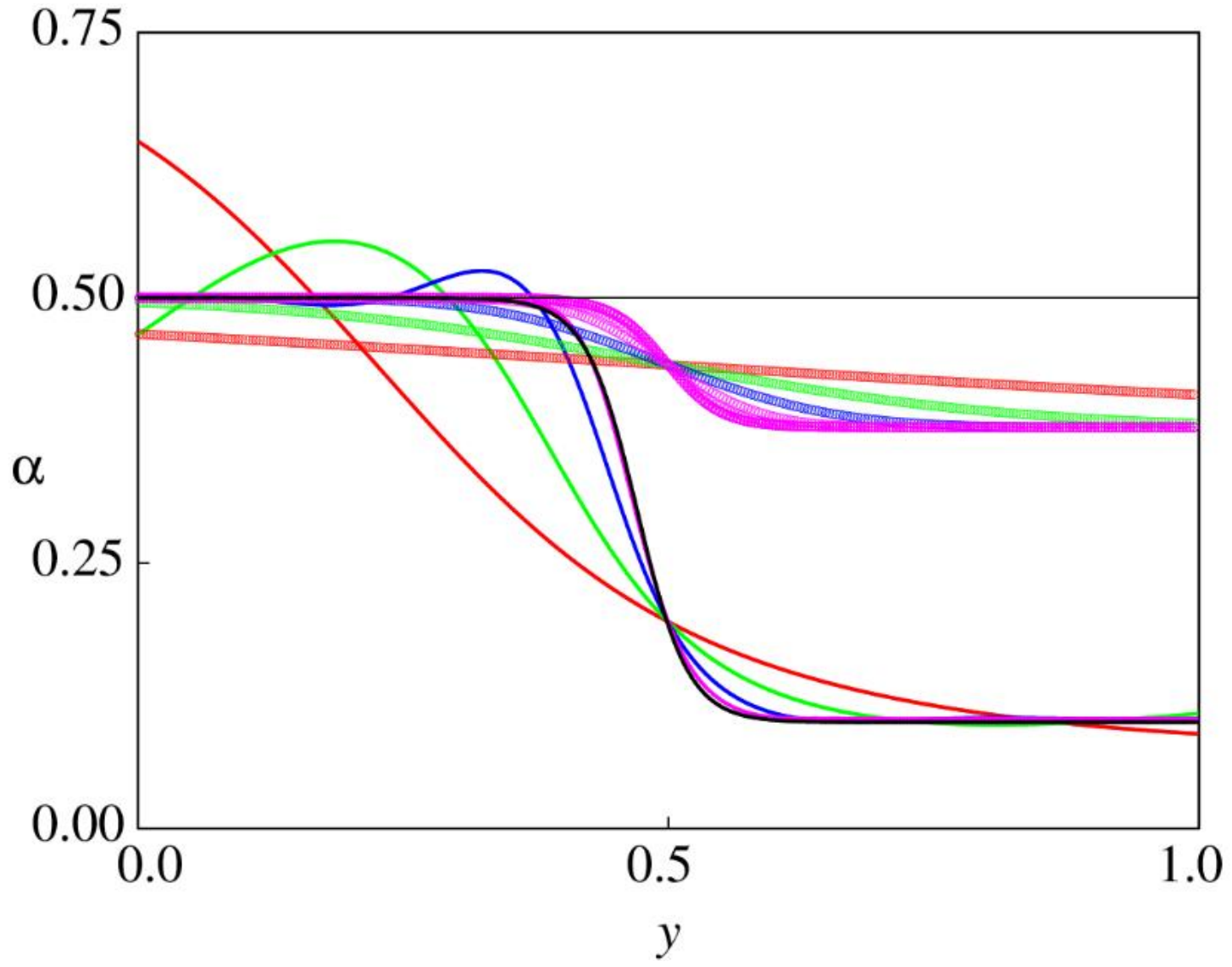

**Figure 10.** Computed average fishing pressure $\alpha = \alpha_y$ in Case (D) for different values of the parameter $\theta$ in the graphon $W$. Black curve represents the average of Nash equilibria. Colored results correspond to: $\theta = 2^{-1}$ (Red), $\theta = 2^{-2}$ (Green), $\theta = 2^{-3}$ (Blue), and $\theta = 2^{-4}$ (Magenta). The non-black lines correspond to the discounted logit dynamic and non-black circles to the G-MFLD.

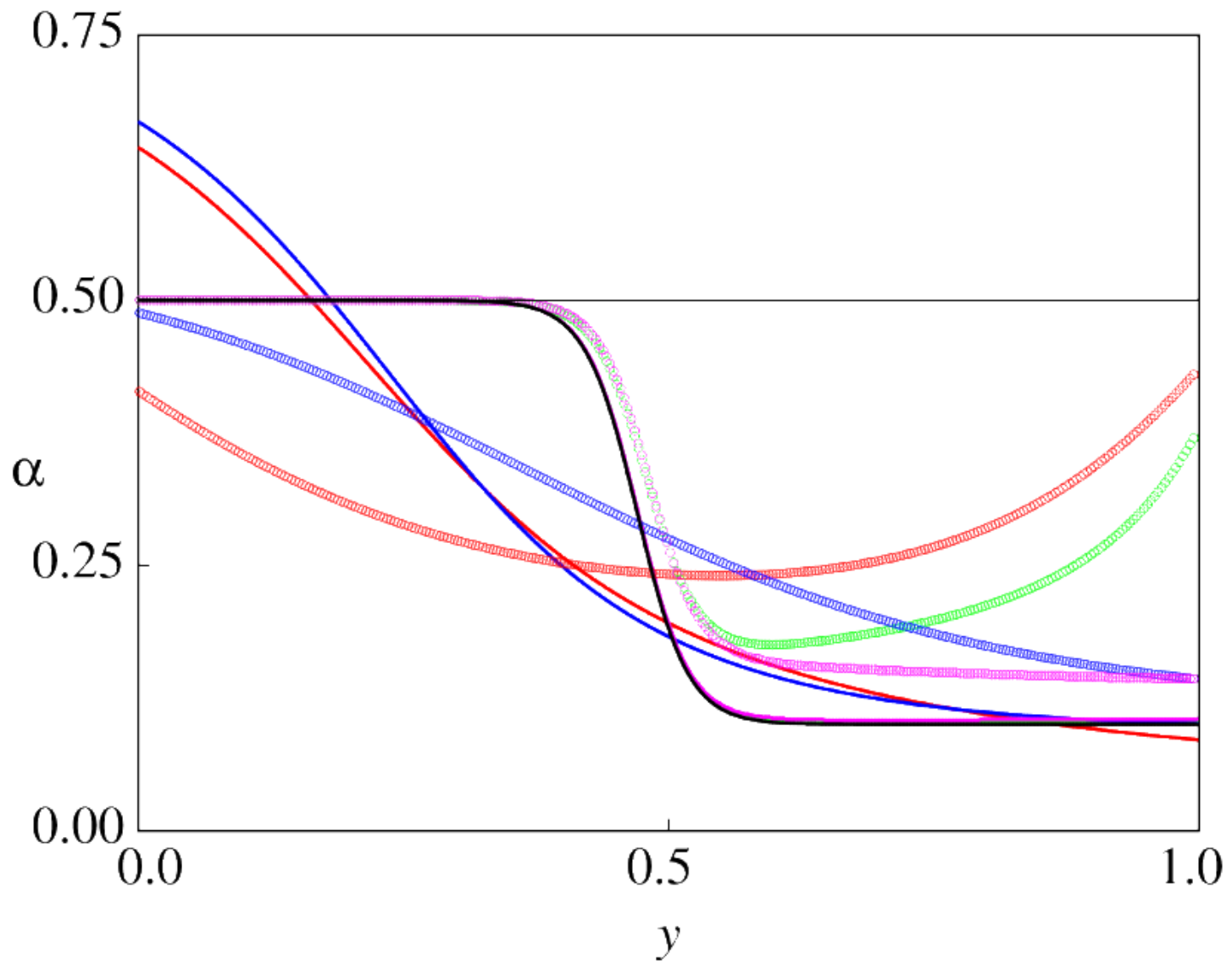

**Figure 11.** Computed average fishing pressure $\alpha = \alpha_y$ for varying $\delta_y$. Colored results correspond to: Case R with $\theta = 2^{-1}$ (Red), Case R with $\theta = 2^{-7}$ (Green), Case M with $\theta = 2^{-1}$ (Blue), and Case M with $\theta = 2^{-7}$ (Magenta). The non-black lines correspond to the discounted logit dynamic and non-black circles to the G-MFLD. Black curve represents the average of Nash equilibria.

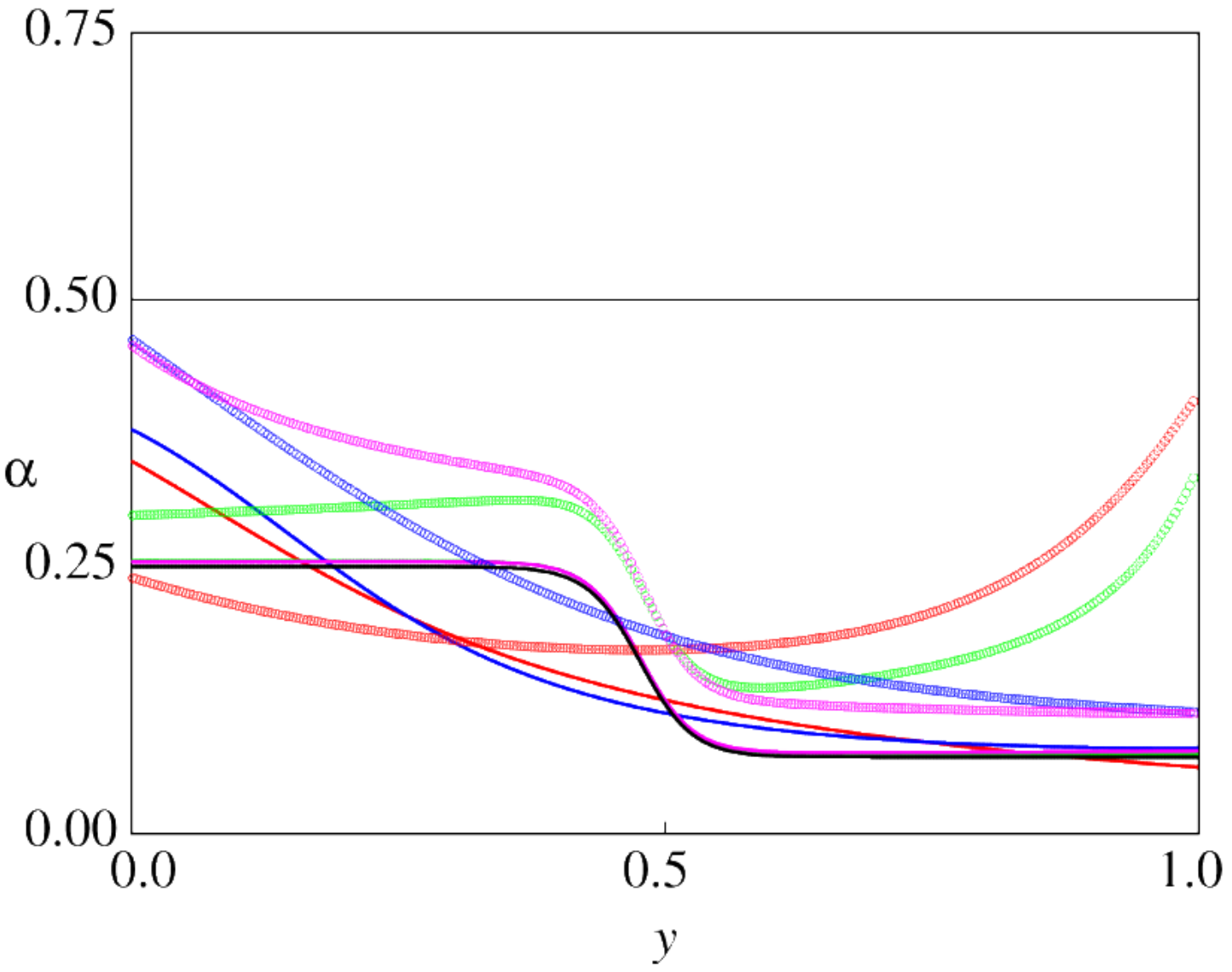

**Figure 12.** Computed average fishing pressure $\alpha = \alpha_y$ for varying $\delta_y$ with larger costs. Figure legends are the same with **Figure 11**.

## 5. Conclusion

The G-MFLD, formulated as a continuum of HJB systems interconnected through a graphon, was derived and analyzed with respect to its well-posedness and regularity. A finite difference scheme for computing the HJB system was presented, along with proofs ensuring the unique existence of bounded solutions and convergence results of the numerical method. Applying the G-MFLD to the case of charr in a mountainous river environment suggested the importance of protecting the habitats of the spotless charr, a rare species, to support sustainable tourism management in the region.

In this study, we considered graphons as bivariate functions; however, an alternative formulation using extended graphons as measures is also possible [88]. Theoretical and numerical analyses would need to be adapted accordingly to address this more complex framework. Additionally, variants such as digraphon, directional versions of graphons [89], and hypergraphons, which capture higher-order interactions [90], can be explored for broader theoretical development. Note that finite difference schemes encounter the curse of dimensionality in high-dimensional or higher-order problems; our scheme is not an exception. Moreover, mathematical properties of the proposed model in such advanced cases may be significantly different from the present low-dimensional case. These theoretical issues are currently open. Finally, a promising application of our graphon mean-field games framework would be modeling a river system as a single entity parameterized through graphon, enabling spatially explicit strategies for environment and fisheries management. We consider that the symmetry assumption of graphon, which was employed for modeling simplicity in this paper, is worth investigating by social surveys in the study area, with which different agents' behavior may emerge.

## Appendix

### A.1 Proofs

***Proof of Proposition 1***

We divide the proof into several steps.

***Step 1: a priori estimate***

We obtain an *a priori* estimate for solutions to the HJB system (23). For any solution $\Phi \in B(\Xi)$ to (23), set $\bar{\Phi} = \sup_{(x,y)\in\Xi} \Phi(x,y)$ and $\underline{\Phi} = \inf_{(x,y)\in\Xi} \Phi(x,y)$. Then, at any point $(x,y)\in\Xi$, we have

$$\begin{aligned}\Phi(x,y) &\le \bar{U} + \frac{1}{\eta_y \delta_y} \ln\left(\int_\Omega e^{\eta_y\{\Phi(z,y)-\Phi(x,y)\}} \mathrm{d}z\right) \\ &\le \bar{U} + \frac{1}{\eta_y \delta_y} \ln\left(\int_\Omega e^{\eta_y \bar{\Phi}} \mathrm{d}z\right) - \frac{1}{\delta_y}\Phi(x,y) \\ &= \bar{U} + \frac{1}{\delta_y}\bar{\Phi} - \frac{1}{\delta_y}\Phi(x,y)\end{aligned} \tag{56}$$

and similarly,

$$\begin{aligned}\Phi(x,y) &\ge 0 + \frac{1}{\eta_y \delta_y} \ln\left(\int_\Omega e^{\eta_y\{\Phi(z,y)-\Phi(x,y)\}} \mathrm{d}z\right) \\ &\ge \frac{1}{\eta_y \delta_y} \ln\left(\int_\Omega e^{\eta_y \underline{\Phi}} \mathrm{d}z\right) - \frac{1}{\delta_y}\Phi(x,y) \\ &\ge \frac{1}{\delta_y}\underline{\Phi} - \frac{1}{\delta_y}\Phi(x,y).\end{aligned} \tag{57}$$

Consequently, for any $(x,y)\in\Xi$, we have

$$\left(1+\frac{1}{\delta_y}\right)\Phi(x,y) \le \bar{U} + \frac{1}{\delta_y}\bar{\Phi} \quad \text{and} \quad \left(1+\frac{1}{\delta_y}\right)\Phi(x,y) \ge \frac{1}{\delta_y}\underline{\Phi}\,, \tag{58}$$

and hence

$$\frac{1}{1+\delta_y}\underline{\Phi} \le \Phi(x,y) \le \frac{1}{1+\delta_y}\left(\delta_y \bar{U} + \bar{\Phi}\right). \tag{59}$$

Taking the infimum on the left-hand side of inequality (59) yields

$$\frac{1}{1+\bar{\delta}}\underline{\Phi} \le \underline{\Phi} \quad \text{and hence} \quad \underline{\Phi} \ge 0\,. \tag{60}$$

Similarly, taking the supremum on the right-hand side of inequality (59) yields

$$\bar{\Phi} \le \frac{1}{1+\underline{\delta}}\left(\bar{\delta}\bar{U} + \bar{\Phi}\right) \quad \text{and hence} \quad \bar{\Phi} \le \frac{\bar{\delta}}{\underline{\delta}}\bar{U}\,. \tag{61}$$

Consequently, by (59)-(61), we obtain the *a priori* bound

$$0 \le \Phi(x,y) \le \frac{\bar{\delta}}{\underline{\delta}}\bar{U}\,. \tag{62}$$

***Step 2: Modified equation***

We consider the following modified equation of the HJB system: for all $(x,y)\in\Xi$ and $\Phi\in B(\Xi)$,

$$\begin{aligned}\Phi(x,y)&=U\left(x,y,\frac{\delta_y}{\delta_y+1}\mu_0(y,\mathrm{d}x)+\frac{1}{\delta_y+1}\frac{1}{\int_\Omega e^{\eta_y(\hat{\Phi}(z,y)-\hat{\Phi}(x,y))}\mathrm{d}z}\mathrm{d}x\right)+\frac{1}{\delta_y\eta_y}\ln\left(\int_\Omega e^{\eta_y\{\hat{\Phi}(z,y)-\hat{\Phi}(x,y)\}}\mathrm{d}z\right)\\ &\equiv\hat{\mathbb{H}}_1[\Phi](x,y)+\hat{\mathbb{H}}_2[\Phi](x,y)\\ &\equiv\hat{\mathbb{H}}[\Phi](x,y),\end{aligned}\tag{63}$$

where $\hat{\Phi}=\max\left\{0,\min\left\{\Phi,\frac{\bar{\delta}}{\underline{\delta}}\bar{U}\right\}\right\}$ is a truncated version of $\Phi$ based on the estimate (62). In the sequel, we use the elementary inequality $|\hat{a}-\hat{b}|\leq|a-b|$ for all $a,b\in\mathbb{R}$.

We fix $(x,y)\in\Xi$. For any $\Phi_1,\Phi_2\in B(\Xi)$, by **Assumption 1**, we have

$$\begin{aligned}\left|\hat{\mathbb{H}}_1[\Phi_1](x,y)-\hat{\mathbb{H}}_1[\Phi_2](x,y)\right|&\leq\left|\begin{aligned}&U\left(x,y,\frac{\delta_y}{\delta_y+1}\mu_0(y,\mathrm{d}x)+\frac{1}{\delta_y+1}\frac{1}{\int_\Omega e^{\eta_y(\hat{\Phi}_1(z,y)-\hat{\Phi}_1(x,y))}\mathrm{d}z}\mathrm{d}x\right)\\ &-U\left(x,y,\frac{\delta_y}{\delta_y+1}\mu_0(y,\mathrm{d}x)+\frac{1}{\delta_y+1}\frac{1}{\int_\Omega e^{\eta_y(\hat{\Phi}_2(z,y)-\hat{\Phi}_2(x,y))}\mathrm{d}z}\mathrm{d}x\right)\end{aligned}\right|\\ &\leq\frac{L_U}{\underline{\delta}+1}\sup_{y\in I}\left\|\frac{e^{\eta_y\hat{\Phi}_1(x,y)}\mathrm{d}x}{\int_\Omega e^{\eta_y\hat{\Phi}_1(z,y)}\mathrm{d}z}-\frac{e^{\eta_y\hat{\Phi}_2(x,y)}\mathrm{d}x}{\int_\Omega e^{\eta_y\hat{\Phi}_2(z,y)}\mathrm{d}z}\right\|_{\mathrm{TV}}.\end{aligned}\tag{64}$$

For all $y\in I$, we have

$$\begin{aligned}&\left\|\frac{e^{\eta_y\hat{\Phi}_1(x,y)}\mathrm{d}x}{\int_\Omega e^{\eta_y\hat{\Phi}_1(z,y)}\mathrm{d}z}-\frac{e^{\eta_y\hat{\Phi}_2(x,y)}\mathrm{d}x\mathrm{d}x}{\int_\Omega e^{\eta_y\hat{\Phi}_2(z,y)}\mathrm{d}z}\right\|_{\mathrm{TV}}\\ &\leq\sup_{|f(\cdot)|\leq1}\int_{x\in\Omega}|f(x)|\left|\frac{e^{\eta_y\hat{\Phi}_1(x,y)}}{\int_{z\in\Omega}e^{\eta_y\hat{\Phi}_1(z,y)}\mathrm{d}z}-\frac{e^{\eta_y\hat{\Phi}_2(x,y)}}{\int_{z\in\Omega}e^{\eta_y\hat{\Phi}_2(z,y)}\mathrm{d}z}\right|\mathrm{d}x\\ &\leq\frac{1}{\int_{z\in\Omega}e^{\eta_y\hat{\Phi}_1(z,y)}\mathrm{d}z\int_{z\in\Omega}e^{\eta_y\hat{\Phi}_2(z,y)}\mathrm{d}z}\int_{x\in\Omega}\left|\int_{z\in\Omega}e^{\eta_y\hat{\Phi}_2(z,y)}\mathrm{d}ze^{\eta_y\hat{\Phi}_1(x,y)}-\int_{z\in\Omega}e^{\eta_y\hat{\Phi}_1(z,y)}\mathrm{d}ze^{\eta_y\hat{\Phi}_2(x,y)}\right|\mathrm{d}x\\ &\leq\int_{x\in\Omega}\left|\int_{z\in\Omega}e^{\eta_y\hat{\Phi}_2(z,y)}\mathrm{d}ze^{\eta_y\hat{\Phi}_1(x,y)}-\int_{z\in\Omega}e^{\eta_y\hat{\Phi}_1(z,y)}\mathrm{d}ze^{\eta_y\hat{\Phi}_2(x,y)}\right|\mathrm{d}x\\ &\leq\int_{x\in\Omega}\left|\int_{z\in\Omega}e^{\eta_y\hat{\Phi}_2(z,y)}\mathrm{d}ze^{\eta_y\hat{\Phi}_1(x,y)}-\int_{z\in\Omega}e^{\eta_y\hat{\Phi}_1(z,y)}\mathrm{d}ze^{\eta_y\hat{\Phi}_1(x,y)}\right|\mathrm{d}x\\ &+\int_{x\in\Omega}\left|\int_{z\in\Omega}e^{\eta_y\hat{\Phi}_1(z,y)}\mathrm{d}ze^{\eta_y\hat{\Phi}_1(x,y)}-\int_{z\in\Omega}e^{\eta_y\hat{\Phi}_1(z,y)}\mathrm{d}ze^{\eta_y\hat{\Phi}_2(x,y)}\right|\mathrm{d}x\\ &=\int_{x\in\Omega}e^{\eta_y\hat{\Phi}_1(x,y)}\mathrm{d}x\int_{z\in\Omega}\left|e^{\eta_y\hat{\Phi}_2(z,y)}-e^{\eta_y\hat{\Phi}_1(z,y)}\right|\mathrm{d}z+\int_{x\in\Omega}\left|e^{\eta_y\hat{\Phi}_1(x,y)}-e^{\eta_y\hat{\Phi}_2(x,y)}\right|\mathrm{d}x\int_{z\in\Omega}e^{\eta_y\hat{\Phi}_1(z,y)}\mathrm{d}z\\ &=2e^{\bar{\eta}\frac{\bar{\delta}}{\underline{\delta}}\bar{U}}\int_{x\in\Omega}\left|e^{\eta_y\hat{\Phi}_1(x,y)}-e^{\eta_y\hat{\Phi}_2(x,y)}\right|\mathrm{d}x.\end{aligned}\tag{65}$$

We also have

$$\int_{x\in\Omega}\left|e^{\eta_y\hat{\Phi}_1(x,y)}-e^{\eta_y\hat{\Phi}_2(x,y)}\right|\mathrm{d}x\le\max_{x\in\Omega}\left|e^{\eta_y\hat{\Phi}_1(x,y)}-e^{\eta_y\hat{\Phi}_2(x,y)}\right|\le\bar{\eta}e^{\bar{\eta}\frac{\bar{\delta}}{\underline{\delta}}\bar{U}}\max_{x\in\Omega}\left|\Phi_1(x,y)-\Phi_2(x,y)\right|. \tag{66}$$

Combining (64)-(66) and using the arbitrariness of $(x,y)\in\Xi$, yields the Lipschitz continuity of $\hat{\mathbb{H}}_1$:

$$\left\|\hat{\mathbb{H}}_1[\Phi_1]-\hat{\mathbb{H}}_1[\Phi_2]\right\|_\infty\le\frac{2\bar{\eta}L_U}{\underline{\delta}+1}e^{2\bar{\eta}\frac{\bar{\delta}}{\underline{\delta}}\bar{U}}\left\|\Phi_1-\Phi_2\right\|_\infty. \tag{67}$$

Similarly, we also obtain the boundedness of $\hat{\mathbb{H}}_1$:

$$\left\|\hat{\mathbb{H}}_1[\Phi_1]\right\|_\infty\le\bar{U}. \tag{68}$$

We next estimate $\hat{\mathbb{H}}_2$. We fix $(x,y)\in\Xi$. For any $\Phi_1,\Phi_2\in B(\Xi)$, by **Assumption 1**, we have

$$\begin{aligned}
&\left|\hat{\mathbb{H}}_2[\Phi_1](x,y)-\hat{\mathbb{H}}_2[\Phi_2](x,y)\right|\\
&=\left|\frac{1}{\delta_y\eta_y}\left\{\ln\left(\int_\Omega e^{\eta_y\hat{\Phi}_1(z,y)}\mathrm{d}z\right)-\ln\left(\int_\Omega e^{\eta_y\hat{\Phi}_2(z,y)}\mathrm{d}z\right)\right\}-\frac{1}{\delta_y}\left(\hat{\Phi}_1(x,y)-\hat{\Phi}_2(x,y)\right)\right|\\
&\le\frac{1}{\underline{\delta}}\frac{1}{\eta_y}\frac{1}{\min\left\{\int_\Omega e^{\eta_y\hat{\Phi}_1(z,y)}\mathrm{d}z,\int_\Omega e^{\eta_y\hat{\Phi}_2(z,y)}\mathrm{d}z\right\}}\left|\int_\Omega e^{\eta_y\hat{\Phi}_1(z,y)}\mathrm{d}z-\int_\Omega e^{\eta_y\hat{\Phi}_2(z,y)}\mathrm{d}z\right|+\frac{1}{\underline{\delta}}\left|\hat{\Phi}_1(x,y)-\hat{\Phi}_2(x,y)\right| \quad (69)\\
&\le\frac{1}{\underline{\delta}}\frac{1}{\eta_y}e^{\bar{\eta}\frac{\bar{\delta}}{\underline{\delta}}\bar{U}}\int_\Omega\left|\eta_y\hat{\Phi}_1(z,y)-\eta_y\hat{\Phi}_2(z,y)\right|\mathrm{d}z+\frac{1}{\underline{\delta}}\left\|\Phi_1-\Phi_2\right\|_\infty\\
&\le\frac{1}{\underline{\delta}}\left(1+e^{\bar{\eta}\frac{\bar{\delta}}{\underline{\delta}}\bar{U}}\right)\left\|\Phi_1-\Phi_2\right\|_\infty.
\end{aligned}$$

Then, by the arbitrariness of $(x,y)\in\Xi$, we obtain the Lipschitz continuity of $\hat{\mathbb{H}}_2$:

$$\left\|\hat{\mathbb{H}}_2[\Phi_1]-\hat{\mathbb{H}}_2[\Phi_2]\right\|_\infty\le\frac{1}{\underline{\delta}}\left(1+e^{\bar{\eta}\frac{\bar{\delta}}{\underline{\delta}}\bar{U}}\right)\left\|\Phi_1-\Phi_2\right\|_\infty. \tag{70}$$

Similarly, we obtain the boundedness of $\hat{\mathbb{H}}_2$:

$$\left\|\hat{\mathbb{H}}_2[\Phi_1]\right\|_\infty\le\sup_{y\in I}\frac{1}{\underline{\delta}\eta_y}\ln\left(\int_\Omega e^{\eta_y\left\{\frac{\bar{\delta}}{\underline{\delta}}\bar{U}-0\right\}}\mathrm{d}z\right)=\frac{\bar{\delta}}{\underline{\delta}^2}\bar{U}. \tag{71}$$

The operator $\hat{\mathbb{H}}$ is bounded due to (68) and (71). Moreover, from (67) and (70), we have

$$\begin{aligned}
\left\|\hat{\mathbb{H}}[\Phi_1]-\hat{\mathbb{H}}[\Phi_2]\right\|_\infty&\le\left\|\hat{\mathbb{H}}_1[\Phi_1]-\hat{\mathbb{H}}_1[\Phi_2]\right\|_\infty+\left\|\hat{\mathbb{H}}_2[\Phi_1]-\hat{\mathbb{H}}_2[\Phi_2]\right\|_\infty\\
&\le\left\{\frac{2\bar{\eta}L_U}{\underline{\delta}+1}e^{2\bar{\eta}\frac{\bar{\delta}}{\underline{\delta}}\bar{U}}+\frac{1}{\underline{\delta}}\left(1+e^{\bar{\eta}\frac{\bar{\delta}}{\underline{\delta}}\bar{U}}\right)\right\}\left\|\Phi_1-\Phi_2\right\|_\infty.
\end{aligned} \tag{72}$$

Therefore, the operator $\hat{\mathbb{H}}$, viewed as a mapping from $B(\Xi)$ to $B(\Xi)$, becomes strictly contractive if

$$\frac{2\bar{\eta}L_U}{\underline{\delta}+1}e^{2\bar{\eta}\frac{\bar{\delta}}{\underline{\delta}}\bar{U}}+\frac{1}{\underline{\delta}}\left(1+e^{\bar{\eta}\frac{\bar{\delta}}{\underline{\delta}}\bar{U}}\right)\in(0,1). \tag{73}$$

This condition corresponds exactly to (25). We then apply Banach's fixed-point theorem (e.g., Theorem 5.7 in Brezis [91]) to the modified equation (63), which ensures the existence of a unique solution in $B(\Xi)$. We denote this solution by $\Psi$ in the sequel.

***Step 3: Uniqueness of the HJB system***

The solution $\Psi$ to the modified equation (63) also solves the original HJB system (23), as it satisfies the *a priori bound* (62) through analytical calculations analogous to those in **Step 1**. This implies that the HJB system (23) admits a unique solution in $B(\Xi)$, namely $\Psi$, that satisfies the bound (62). Indeed, if there were more than one such solutions, say $\Psi_1, \Psi_2$, then we have

$$\left\|\Psi_1 - \Psi_2\right\|_\infty = \left\|\mathbb{H}[\Psi_1] - \mathbb{H}[\Psi_2]\right\|_\infty \leq \left\{ \frac{2\bar{\eta} L_U}{\underline{\delta}+1} e^{2\bar{\eta}\frac{\bar{\delta}}{\underline{\delta}}\bar{U}} + \frac{1}{\underline{\delta}}\left(1 + e^{\bar{\eta}\frac{\bar{\delta}}{\underline{\delta}}\bar{U}}\right)\right\} \left\|\Psi_1 - \Psi_2\right\|_\infty < \left\|\Psi_1 - \Psi_2\right\|_\infty, \quad (74)$$

leading to a contradiction, proving the uniqueness. The proof is then completed.

□

***Proof of Proposition 2***

We set $\Delta = \left|\Phi(x_1, y) - \Phi(x_2, y)\right|$. By (23) and **Assumption 2**, we have

$$\begin{aligned} \Delta &\leq \left|\mathbb{H}_1[\Phi](x_1, y) - \mathbb{H}_1[\Phi](x_2, y)\right| + \left|\mathbb{H}_2[\Phi](x_1, y) - \mathbb{H}_2[\Phi](x_2, y)\right| \\ &\leq \phi_U\left(\left|x_1 - x_2\right|\right) + \left|\mathbb{H}_2[\Phi](x_1, y) - \mathbb{H}_2[\Phi](x_2, y)\right|. \end{aligned} \quad (75)$$

The second term can be rewritten as

$$\begin{aligned} \left|\mathbb{H}_2[\Phi](x_1, y) - \mathbb{H}_2[\Phi](x_2, y)\right| &= \frac{1}{\delta_y \eta_y}\left|\ln\left(\int_\Omega e^{\eta_y \Phi(z,y)} \mathrm{d}z\right) - \ln\left(\int_\Omega e^{\eta_y \Phi(z,y)} \mathrm{d}z\right) - \eta_y \Phi(x_1, y) + \eta_y \Phi(x_2, y)\right| \\ &= \frac{1}{\delta_y}\Delta. \end{aligned} \quad (76)$$

Combining (75) and (76) yields

$$\Delta \leq \phi_U\left(\left|x_1 - x_2\right|\right) + \frac{1}{\delta_y}\Delta, \quad (77)$$

and hence the desired result (26) follows from $\delta_y \geq \underline{\delta} > 1$ for all $y \in I$.

□

***Proof of Proposition 3***

We set $\Delta = \left|\Phi(x, y_1) - \Phi(x, y_2)\right|$. By (23) and **Assumption 3**, we have

$$\begin{aligned} \Delta &\leq \left|\mathbb{H}_1[\Phi](x, y_1) - \mathbb{H}_1[\Phi](x, y_2)\right| + \left|\mathbb{H}_2[\Phi](x, y_1) - \mathbb{H}_2[\Phi](x, y_2)\right| \\ &\leq \psi_U\left(\left|y_1 - y_2\right|\right) + \left|\mathbb{H}_2[\Phi](x, y_1) - \mathbb{H}_2[\Phi](x, y_2)\right|. \end{aligned} \quad (78)$$

The last term is evaluated as

$$\begin{aligned}\left|\mathbb{H}_2[\Phi](x,y_1)-\mathbb{H}_2[\Phi](x,y_2)\right| &\le \frac{1}{\delta\eta}\left|\ln\left(\int_\Omega e^{\eta\Phi(z,y_1)}\mathrm{d}z\right)-\ln\left(\int_\Omega e^{\eta\Phi(z,y_2)}\mathrm{d}z\right)\right|+\frac{1}{\delta}\Delta \\ &\le \frac{1}{\delta\eta}\left|\int_\Omega e^{\eta\Phi(z,y_1)}\mathrm{d}z-\int_\Omega e^{\eta\Phi(z,y_2)}\mathrm{d}z\right|+\frac{1}{\delta}\Delta \\ &\le \frac{1}{\delta\eta}\int_\Omega\left|e^{\eta\Phi(z,y_1)}-e^{\eta\Phi(z,y_2)}\right|\mathrm{d}z+\frac{1}{\delta}\Delta \\ &\le \frac{1}{\delta}e^{\eta\bar{U}}\int_\Omega\left|\Phi(z,y_1)-\Phi(z,y_2)\right|\mathrm{d}z+\frac{1}{\delta}\Delta.\end{aligned} \tag{79}$$

We set $\bar{\Delta}=\sup_{x\in\Omega}\left|\Phi(x,y_1)-\Phi(x,y_2)\right|$. By (79), we obtain

$$\left|\mathbb{H}_2[\Phi](x,y_1)-\mathbb{H}_2[\Phi](x,y_2)\right|\le\frac{1}{\delta}e^{\eta\bar{U}}\bar{\Delta}+\frac{1}{\delta}\Delta. \tag{80}$$

Combining (78) and (80) yields

$$\left(1-\frac{1}{\delta}\right)\Delta\le\psi_U\left(|y_1-y_2|\right)+\frac{1}{\delta}e^{\eta\bar{U}}\bar{\Delta}. \tag{81}$$

Taking the supremum over $x\in\Omega$ on both sides of (81) and rearranging yields

$$\left(1-\frac{1}{\delta}\left(1+e^{\eta\bar{U}}\right)\right)\bar{\Delta}\le\psi_U\left(|y_1-y_2|\right), \tag{82}$$

which implies the desired result (27), where $\frac{1}{\delta}\left(1+e^{\eta\bar{U}}\right)\in(0,1)$ due to (25).

□

***Proof of Proposition 5***

Fix $(x,y)\in\Xi$. By **Proposition 4**, the solution $\Phi$ to the HJB system is equi-continuous with respect to $\delta$ provided is sufficiently large. Indeed, we have $\delta/(\delta-1)\le 2$ and $1-\frac{1}{\delta}\left(1+e^{\eta\bar{U}}\right)\ge\frac{1}{2}$ for sufficiently large $\delta>0$, which ensures equi-continuity for such $\delta$:

$$\left|\Phi(x_1,y_1)-\Phi(x_2,y_2)\right|\le 2\left(\phi_U\left(|x_1-x_2|\right)+\psi_U\left(|y_1-y_2|\right)\right). \tag{83}$$

Moreover, $\Phi$ is uniformly bounded for any sufficiently large $\delta$. Then, by the classical Ascoli-Arzela theorem, there exists a subsequence $\delta_k\to+\infty$ ($k=1,2,3,...$) of $\delta$ such that there exists a function $P=\lim_{k\to+\infty}\Phi_k$ (we write $\Phi_k=\Phi_{\delta=\delta_k}$) in $\Xi$ with the estimate (83). For each $(x,y)\in\Xi$, we obtain

$$\lim_{k\to+\infty}\frac{1}{\delta_k\eta_y}\ln\left(\int_\Omega e^{\eta_y\{\Phi_k(z,y)-\Phi_k(x,y)\}}\mathrm{d}z\right)=0 \tag{84}$$

and

$$
\begin{aligned}
&\lim_{k\to+\infty} U\left(x, y, \frac{\delta_k}{\delta_k+1}\mu_0(y,\mathrm{d}x) + \frac{1}{\delta_k+1}\frac{1}{\int_\Omega e^{\eta_y(\Phi_k(z,y)-\Phi_k(x,y))}\mathrm{d}z}\mathrm{d}x\right) \\
&= U\left(x, y, \lim_{k\to+\infty}\frac{\delta_k}{\delta_k+1}\mu_0(y,\mathrm{d}x) + \lim_{k\to+\infty}\frac{1}{\delta_k+1}\frac{1}{\int_\Omega e^{\eta_y(\Phi_k(z,y)-\Phi_k(x,y))}\mathrm{d}z}\mathrm{d}x\right) \\
&= U(x, y, \mu_0).
\end{aligned}
\tag{85}
$$

Moreover, by applying the classical dominated convergence theorem, we obtain the estimate

$$
\begin{aligned}
&\lim_{k\to+\infty}\left|U(x,y,m) - U(x,y,\mu_0)\right| \\
&\le \lim_{k\to+\infty} L_U \left\|m-\mu_0\right\|_{\mathrm{TV}}^{(I)} \\
&= L_U \lim_{\delta\to+\infty}\left\|\frac{\delta_k}{\delta_k+1}\mu_0(y,\mathrm{d}x) + \frac{1}{\delta_k+1}\frac{1}{\int_\Omega e^{\eta_y(\Phi_k(z,y)-\Phi_k(x,y))}\mathrm{d}z}\mathrm{d}x - \mu_0(y,\mathrm{d}x)\right\|_{\mathrm{TV}}^{(I)} \\
&= 0.
\end{aligned}
\tag{86}
$$

Combining (84)-(86) yields $P(x,y) = U(x,y,\mu_0)$ at each $(x,y)\in\Xi$, and its right-hand side satisfies (83) by **Assumptions 2 and 3**. Due to (84)-(86) and the uniform boundedness of $\Phi$ for large $\delta > 0$, these results hold for any subsequences $\delta_k \to +\infty$ ($k = 1,2,3,...$), hence concluding the proof.

□

***Proof of Proposition 6***

Assume that $\Phi_{i,j}$ is minimized at some $(i,j) = (i',j')$. Then, we have

$$
\begin{aligned}
\Phi_{i',j'} &\ge 0 + \frac{1}{\delta_{j'}\eta_{j'}}\ln\left(\sum_{k=1}^{N_x} e^{\eta_{j'}\{\Phi_{k,j'}-\Phi_{i',j'}\}}\Delta x\right) \\
&\ge \frac{1}{\delta_{j'}\eta_{j'}}\ln\left(\sum_{k=1}^{N_x} e^{\eta_{j'}\{\Phi_{i',j'}-\Phi_{i',j'}\}}\Delta x\right) \\
&= \frac{1}{\delta_{j'}\eta_{j'}}\ln\left(\sum_{k=1}^{N_x}\Delta x\right) \\
&= 0.
\end{aligned}
\tag{87}
$$

Similarly, assume that $\Phi_{i,j}$ is maximized at some $(i,j) = (i',j')$. Then, we have

$$
\begin{aligned}
\Phi_{i',j'} &\le \bar{U} + \frac{1}{\delta_{j'}\eta_{j'}}\ln\left(\sum_{k=1}^{N_x} e^{\eta_{j'}\{\Phi_{k,j'}-\Phi_{i',j'}\}}\Delta x\right) \\
&\le \bar{U} + \frac{1}{\delta_{j'}\eta_{j'}}\ln\left(\sum_{k=1}^{N_x} e^{\eta_{j'}\{\Phi_{i',j'}-\Phi_{i',j'}\}}\Delta x\right) \\
&= \bar{U}.
\end{aligned}
\tag{88}
$$

The proof is completed by combining (87) and (88).

□

***Proof of Proposition 8***

First, we recall the original and discretized HJB systems:

**(HJB system)**

$$\Phi(x,y)=\int_{w\in I} g\left(x,w,\frac{1}{\delta+1}\frac{\int_{q\in\Omega} h(q)e^{\eta\Phi(q,w)}\mathrm{d}q}{\int_{z\in\Omega} e^{\eta\Phi(z,w)}\mathrm{d}z}\right)W(w,y)\mathrm{d}w+\frac{1}{\delta\eta}\ln\left(\int_{\Omega} e^{\eta\{\Phi(z,y)-\Phi(x,y)\}}\mathrm{d}z\right) \tag{89}$$

for all $x,y\in\Xi$. The solution is continuous (under **Assumptions 2 and 3**) and bounded between 0 and $\bar{U}$.

**(Discretized HJB system)**

$$\Phi_{i,j}=\sum_{l=1}^{N} g\left(x_i,y_l,\frac{1}{\delta+1}\frac{\sum_{k=1}^{N} h(x_k)e^{\eta\Phi_{k,l}}\Delta x}{\sum_{k=1}^{N_x} e^{\eta\Phi_{k,l}}\Delta x}\right)W_{l,j}\Delta y+\frac{1}{\delta\eta}\ln\left(\sum_{k=1}^{N} e^{\eta\{\Phi_{k,j}-\Phi_{i,j}\}}\Delta x\right) \tag{90}$$

at each $P_{i,j}$ ($1\le i,j\le N$). The solution is bounded between 0 and $\bar{U}$.

The solution obtained from (90) is extended to a piecewise-constant function $\Psi$ on $B(\Xi)$ such that $\Psi(x,y)=\Phi_{i,j}$ if $(x,y)\in C_{i,j}\equiv\Omega_i\times\Omega_j$. Similarly, we define a piecewise-constant function such that $\hat{h}(\cdot)=h(x_i)$ if $(x,y)\in C_{i,j}$. We also define a piecewise-constant function such that $\hat{g}(x,y,\cdot)=g(x_i,y_j,\cdot)$ if $(x,y)\in C_{i,j}$.

Under these settings, if $(x,y)\in C_{i,j}$, then

$$\frac{\sum_{k=1}^{N} h(x_k)e^{\eta\Phi_{k,l}}\Delta x}{\sum_{k=1}^{N} e^{\eta\Phi_{k,l}}\Delta x}=\frac{\int_{z\in\Omega}\hat{h}(z)e^{\eta\Psi(z,y_l)}\mathrm{d}z}{\int_{z\in\Omega} e^{\eta\Psi(z,y_l)}\mathrm{d}z} \tag{91}$$

and

$$\frac{1}{\delta\eta}\ln\left(\sum_{k=1}^{N} e^{\eta\{\Phi_{k,j}-\Phi_{i,j}\}}\Delta x\right)=\frac{1}{\delta\eta}\ln\left(\int_{z\in\Omega} e^{\eta\{\Psi(z,y)-\Psi(x,y)\}}\mathrm{d}z\right). \tag{92}$$

We also have

$$
\begin{aligned}
&\sum_{l=1}^{N}\hat{g}\left(x_i,y_l,\frac{1}{\delta+1}\frac{\sum_{k=1}^{N}h(x_k)e^{\eta\Phi_{k,l}}\Delta x}{\sum_{k=1}^{N}e^{\eta\Phi_{k,l}}\Delta x}\right)W_{l,j}\Delta y\\
&=\sum_{l=1}^{N}\hat{g}\left(x_i,y_l,\frac{1}{\delta+1}\frac{\int_{z\in\Omega}\hat{h}(z)e^{\eta\Psi(z,y_l)}\mathrm{d}z}{\int_{z\in\Omega}e^{\eta\Psi(z,y_l)}\mathrm{d}z}\right)\frac{1}{\Delta y}\left(\int_{y_l-\Delta y/2}^{y_l+\Delta y/2}W(w,y_j)\mathrm{d}w\right)\Delta y\\
&=\sum_{l=1}^{N}\hat{g}\left(x_i,y_l,\frac{1}{\delta+1}\frac{\int_{z\in\Omega}\hat{h}(z)e^{\eta\Psi(z,y_l)}\mathrm{d}z}{\int_{z\in\Omega}e^{\eta\Psi(z,y_l)}\mathrm{d}z}\right)\int_{y_l-\Delta y/2}^{y_l+\Delta y/2}W(w,y_j)\mathrm{d}w \qquad (93)\\
&=\sum_{l=1}^{N}\int_{y_l-\Delta y/2}^{y_l+\Delta y/2}\hat{g}\left(x_i,y_l,\frac{1}{\delta+1}\frac{\int_{z\in\Omega}\hat{h}(z)e^{\eta\Psi(z,y_l)}\mathrm{d}z}{\int_{z\in\Omega}e^{\eta\Psi(z,y_l)}\mathrm{d}z}\right)W(w,y_j)\mathrm{d}w\\
&=\int_{w\in I}\hat{g}\left(x_i,w,\frac{1}{\delta+1}\frac{\int_{z\in\Omega}\hat{h}(z)e^{\eta\Psi(z,w)}\mathrm{d}z}{\int_{z\in\Omega}e^{\eta\Psi(z,w)}\mathrm{d}z}\right)W(w,y_j)\mathrm{d}w.
\end{aligned}
$$

Therefore, (90) can be rewritten as

$$
\Psi(x,y)=\int_{w\in I}\hat{g}\left(x_i,w,\frac{1}{\delta+1}\frac{\int_{z\in\Omega}\hat{h}(z)e^{\eta\Psi(z,w)}\mathrm{d}z}{\int_{z\in\Omega}e^{\eta\Psi(z,w)}\mathrm{d}z}\right)W(w,y_j)\mathrm{d}w+\frac{1}{\delta\eta}\ln\left(\int_{z\in\Omega}e^{\eta\{\Psi(z,y)-\Psi(x,y)\}}\mathrm{d}z\right) \qquad (94)
$$

for all $(x,y)\in C_{i,j}$ $(1\le i,j\le N)$.

We evaluate the difference between the true and numerical solutions in each $C_{i,j}$. At each $(x,y)\in C_{i,j}$, by (89)-(90), we have

$$
\begin{aligned}
|\Phi-\Psi|(x,y)\le&\left|\begin{array}{l}\int_{w\in I}g\left(x,w,\frac{1}{\delta+1}\frac{\int_{q\in\Omega}h(q)e^{\eta\Phi(q,w)}\mathrm{d}q}{\int_{z\in\Omega}e^{\eta\Phi(z,w)}\mathrm{d}z}\right)W(w,y)\mathrm{d}w\\ -\int_{w\in I}\hat{g}\left(x,w,\frac{1}{\delta+1}\frac{\int_{q\in\Omega}\hat{h}(q)e^{\eta\Psi(q,w)}\mathrm{d}q}{\int_{z\in\Omega}e^{\eta\Psi(z,w)}\mathrm{d}z}\right)W(w,y)\mathrm{d}w\end{array}\right|\\
&+\frac{1}{\delta\eta}\left|\ln\left(\int_{\Omega}e^{\eta\{\Phi(z,y)-\Phi(x,y)\}}\mathrm{d}z\right)-\ln\left(\int_{z\in\Omega}e^{\eta\{\Psi(z,y)-\Psi(x,y)\}}\mathrm{d}z\right)\right| \qquad (95)\\
&=K_1+K_2.
\end{aligned}
$$

To simplify notations, we set

$$
G(w)=\frac{1}{1+\delta}\frac{\int_{q\in\Omega}h(q)e^{\eta\Phi(q,w)}\mathrm{d}q}{\int_{z\in\Omega}e^{\eta\Phi(z,w)}\mathrm{d}z}\quad\text{and}\quad\hat{G}(w)=\frac{1}{\delta+1}\frac{\int_{q\in\Omega}\hat{h}(q)e^{\eta\Psi(q,w)}\mathrm{d}q}{\int_{z\in\Omega}e^{\eta\Psi(z,w)}\mathrm{d}z}. \qquad (96)
$$

For $K_2$, we have

$$
\begin{aligned}
K_2 &= \frac{1}{\delta\eta}\left|\ln\left(\int_\Omega e^{\eta\{\Phi(z,y)-\Phi(x,y)\}}\mathrm{d}z\right)-\ln\left(\int_{z\in\Omega} e^{\eta\{\Psi(z,y)-\Psi(x,y)\}}\mathrm{d}z\right)\right| \\
&\le \frac{1}{\delta\eta}\left|\ln\left(\int_\Omega e^{\eta\Phi(z,y)}\mathrm{d}z\right)-\ln\left(\int_{z\in\Omega} e^{\eta\Psi(z,y)}\mathrm{d}z\right)\right|+\frac{1}{\delta}\left|\Phi-\Psi\right|(x,y) \\
&\le \frac{1}{\delta\eta}\int_\Omega\left|e^{\eta\Phi(z,y)}-e^{\eta\Psi(z,y)}\right|\mathrm{d}z+\frac{1}{\delta}\left|\Phi-\Psi\right|(x,y) \\
&\le \frac{1}{\delta\eta}e^{\eta\bar{U}}\int_\Omega \eta\left|\Phi(z,y)-\Psi(z,y)\right|\mathrm{d}z+\frac{1}{\delta}\left|\Phi-\Psi\right|(x,y) \\
&\le \frac{1}{\delta}\left(1+e^{\eta\bar{U}}\right)\left\|\Phi-\Psi\right\|_\infty .
\end{aligned} \tag{97}
$$

For $K_1$, we have

$$
\begin{aligned}
K_2 &= \left|\int_{w\in I} g(x,w,G(w))W(w,y)\mathrm{d}w-\int_{w\in I}\hat{g}\left(x,w,\hat{G}(w)\right)W\left(w,y_j\right)\mathrm{d}w\right| \\
&\le \left|\int_{w\in I} g(x,w,G(w))W(w,y)\mathrm{d}w-\int_{w\in I} g(x,w,G(w))W\left(w,y_j\right)\mathrm{d}w\right| \\
&\quad +\left|\int_{w\in I} g(x,w,G(w))W\left(w,y_j\right)\mathrm{d}w-\int_{w\in I}\hat{g}\left(x,w,\hat{G}(w)\right)W\left(w,y_j\right)\mathrm{d}w\right| \\
&= K_3+K_4 .
\end{aligned} \tag{98}
$$

We have the estimates

$$
\begin{aligned}
K_3 &\le \int_{w\in I} g(x,w,G(w))\left|W\left(w,y_j\right)-W(w,y)\right|\mathrm{d}w \\
&\le \bar{U}\int_{w\in I}\left|W\left(w,y_j\right)-W(w,y)\right|\mathrm{d}w \\
&\le \bar{U}\sup_j \sup_{y\in\left(y_j-\Delta y/2,\, y_j+\Delta y/2\right)}\int_{w\in I}\left|W\left(w,y_j\right)-W(w,y)\right|\mathrm{d}w
\end{aligned} \tag{99}
$$

and

$$
\begin{aligned}
K_4 &\le \left|\int_{w\in I} g(x,w,G(w))W\left(w,y_j\right)\mathrm{d}w-\int_{w\in I}\hat{g}(x,w,G(w))W\left(w,y_j\right)\mathrm{d}w\right| \\
&\quad +\left|\int_{w\in I}\hat{g}(x,w,G(w))W\left(w,y_j\right)\mathrm{d}w-\int_{w\in I}\hat{g}\left(x,w,\hat{G}(w)\right)W\left(w,y_j\right)\mathrm{d}w\right| \\
&= K_{4,1}+K_{4,2} .
\end{aligned} \tag{100}
$$

For $K_{4,1}$, we have

$$
K_{4,1} \le \int_{w\in I}\left|g(x,w,G(w))-\hat{g}(x,w,G(w))\right|W\left(w,y_j\right)\mathrm{d}w \le \bar{W}\sup_{\substack{x\in\Omega\\ w\in I}}\left|g(x,w,G(w))-\hat{g}(x,w,G(w))\right| . \tag{101}
$$

The last term converges to 0 as $N\to+\infty$. Indeed, due to the Lipschitz continuity of $g$ and the fact that $\hat{g}$ is its piece-wise continuous approximation, we have

$$
\begin{aligned}
\sup_{\substack{x\in\Omega\\ w\in I}}\left|g(x,w,G(w))-\hat{g}(x,w,G(w))\right| &= \sup_{i\le i,j\le N}\ \sup_{(x,y)\in C_{i,j}}\left|g(x,y,G(y))-\hat{g}(x,y,G(y))\right| \\
&= \sup_{i\le i,j\le N} C_g N^{-1}\to 0 \text{ as } N\to+\infty
\end{aligned} \tag{102}
$$

with some constant $C_g>0$ depending only on $g$. For $K_{4,2}$, by the definition of $\hat{g}$ along with its Lipschitz continuity in the third argument, we have

$$
\begin{aligned}
K_{4,2} &\le \int_{w\in I}\left|\hat{g}\left(x,w,G(w)\right)-\hat{g}\left(x,w,\hat{G}(w)\right)\right|W\left(w,y_j\right)\mathrm{d}w \\
&\le L_g\int_{w\in I}\left|G(w)-\hat{G}(w)\right|W\left(w,y_j\right)\mathrm{d}w \\
&\le \frac{L_g\bar{W}}{\delta+1}\sup_{w\in I}\left|\frac{\int_{q\in\Omega}h(q)e^{\eta\Phi(q,w)}\mathrm{d}q}{\int_{z\in\Omega}e^{\eta\Phi(z,w)}\mathrm{d}z}-\frac{\int_{q\in\Omega}\hat{h}(q)e^{\eta\Psi(q,w)}\mathrm{d}q}{\int_{z\in\Omega}e^{\eta\Psi(z,w)}\mathrm{d}z}\right|.
\end{aligned}
\tag{103}
$$

For each $w\in I$, we have

$$
\begin{aligned}
&\left|\frac{\int_{q\in\Omega}h(q)e^{\eta\Phi(q,w)}\mathrm{d}q}{\int_{z\in\Omega}e^{\eta\Phi(z,w)}\mathrm{d}z}-\frac{\int_{q\in\Omega}\hat{h}(q)e^{\eta\Psi(q,w)}\mathrm{d}q}{\int_{z\in\Omega}e^{\eta\Psi(z,w)}\mathrm{d}z}\right| \\
&=\frac{\left|\int_{z\in\Omega}e^{\eta\Psi(z,w)}\mathrm{d}z\int_{q\in\Omega}h(q)e^{\eta\Phi(q,w)}\mathrm{d}q-\int_{z\in\Omega}e^{\eta\Phi(z,w)}\mathrm{d}z\int_{q\in\Omega}\hat{h}(q)e^{\eta\Psi(q,w)}\mathrm{d}q\right|}{\int_{z\in\Omega}e^{\eta\Psi(z,w)}\mathrm{d}z\int_{z\in\Omega}e^{\eta\Phi(z,w)}\mathrm{d}z} \\
&\le\left|\int_{z\in\Omega}e^{\eta\Psi(z,w)}\mathrm{d}z\int_{q\in\Omega}h(q)e^{\eta\Phi(q,w)}\mathrm{d}q-\int_{z\in\Omega}e^{\eta\Psi(z,w)}\mathrm{d}z\int_{q\in\Omega}\hat{h}(q)e^{\eta\Psi(q,w)}\mathrm{d}q\right| \\
&+\left|\int_{z\in\Omega}e^{\eta\Psi(z,w)}\mathrm{d}z\int_{q\in\Omega}\hat{h}(q)e^{\eta\Psi(q,w)}\mathrm{d}q-\int_{z\in\Omega}e^{\eta\Phi(z,w)}\mathrm{d}z\int_{q\in\Omega}\hat{h}(q)e^{\eta\Psi(q,w)}\mathrm{d}q\right| \\
&\le\int_{z\in\Omega}e^{\eta\Psi(z,y)}\mathrm{d}z\left(\int_{q\in\Omega}\left|h(q)e^{\eta\Phi(q,w)}-\hat{h}(q)e^{\eta\Psi(q,w)}\right|\mathrm{d}q\right) \\
&+\int_{q\in\Omega}\hat{h}(q)e^{\eta\Psi(q,y)}\mathrm{d}q\left(\int_{z\in\Omega}\left|e^{\eta\Phi(z,w)}-e^{\eta\Psi(z,w)}\right|\mathrm{d}z\right) \\
&\le\int_{z\in\Omega}e^{\eta\Psi(z,y)}\mathrm{d}z\left(\int_{q\in\Omega}\left|h(q)e^{\eta\Phi(q,w)}-\hat{h}(q)e^{\eta\Phi(q,w)}\right|\mathrm{d}q\right) \\
&+\int_{z\in\Omega}e^{\eta\Psi(z,y)}\mathrm{d}z\left(\int_{q\in\Omega}\left|\hat{h}(q)e^{\eta\Phi(q,w)}-\hat{h}(q)e^{\eta\Psi(q,w)}\right|\mathrm{d}q\right) \\
&+\int_{q\in\Omega}\hat{h}(q)e^{\eta\Psi(q,y)}\mathrm{d}q\left(\int_{z\in\Omega}\left|e^{\eta\Phi(z,w)}-e^{\eta\Psi(z,w)}\right|\mathrm{d}z\right) \\
&\le e^{2\eta\bar{U}}\int_{z\in\Omega}\left|h(z)-\hat{h}(z)\right|\mathrm{d}z+2\bar{h}e^{\eta\bar{U}}\int_{z\in\Omega}\left|e^{\eta\Phi(z,w)}-e^{\eta\Psi(z,w)}\right|\mathrm{d}z \\
&\le e^{2\eta\bar{U}}\int_{z\in\Omega}\left|h(z)-\hat{h}(z)\right|\mathrm{d}z+2\bar{h}\eta e^{2\eta\bar{U}}\left\|\Phi-\Psi\right\|_\infty .
\end{aligned}
\tag{104}
$$

The last line of (104) is independent of $w\in I$.

Consequently, we obtain the inequality

$$
\begin{aligned}
\left|\Phi-\Psi\right|(x,y)&\le\frac{1}{\delta}\left(1+e^{\eta\bar{U}}\right)\left\|\Phi-\Psi\right\|_\infty+\bar{U}\sup_{i,j}\sup_{y\in C_{i,j}}\int_{w\in I}\left|W\left(w,y_j\right)-W\left(w,y\right)\right|\mathrm{d}w \\
&+\bar{W}\sup_{\substack{x\in\Omega\\ w\in I}}\left|g\left(x,w,G(w)\right)-\hat{g}\left(x,w,G(w)\right)\right| \\
&+\frac{L_g\bar{W}}{\delta+1}e^{2\eta\bar{U}}\int_{z\in\Omega}\left|h(z)-\hat{h}(z)\right|\mathrm{d}z+2\frac{L_g\bar{W}}{\delta+1}\bar{h}\eta e^{2\eta\bar{U}}\left\|\Phi-\Psi\right\|_\infty .
\end{aligned}
\tag{105}
$$

By the arbitrariness of $(x,y)$, for a sufficiently large $\delta>0$, we obtain the estimate

$$\left(1-\frac{1}{\delta}\left(1+e^{\eta\bar{U}}\right)-2\frac{L_g\bar{W}}{\delta+1}\bar{h}\eta e^{2\eta\bar{U}}\right)\|\Phi-\Psi\|_\infty \leq \underbrace{\bar{U}\sup_{i,j}\sup_{y\in C_{i,j}}\int_{w\in I}\left|W\left(w,y_j\right)-W\left(w,y\right)\right|\mathrm{d}w}_{\text{Error due to graphon } W}$$

$$+\underbrace{\bar{W}\sup_{\substack{x\in\Omega\\ w\in I}}\left|g\left(x,w,G\left(w\right)\right)-\hat{g}\left(x,w,G\left(w\right)\right)\right|}_{\text{Error due to coefficient } g} \tag{106}$$

$$+\underbrace{\frac{L_g\bar{W}}{\delta+1}e^{2\eta\bar{U}}\int_{z\in\Omega}\left|h\left(z\right)-\hat{h}\left(z\right)\right|\mathrm{d}z}_{\text{Error due to coefficient } h}.$$

The right-hand side of (106) converges to 0 as $N\to+\infty$. Therefore, the conclusion of the proposition due to

$$\sup_{1\leq i,j\leq N}\left|\Phi\left(x_i,y_j\right)-\Phi_{i,j}\right|\leq\sup_{(x,y)\in\Xi}\left|\Phi\left(x,y\right)-\Psi\left(x,y\right)\right|=\|\Phi-\Psi\|_\infty. \tag{107}$$

□

***Proof of Proposition 9***

First, fix $(a,j)$ with $1\leq a\leq N_x$ and $1\leq j\leq N_y$. Concerning the second term of (49), we have

$$\frac{\partial}{\partial\Phi_{a,j}}\left\{-\frac{1}{\delta_j\eta_j}\ln\left(\sum_{k=1}^{N_x}e^{\eta_j\Phi_{k,j}}\Delta x\right)\right\}=-\frac{1}{\delta_j}\frac{e^{\eta_j\Phi_{a,j}}}{\sum_{k=1}^{N_x}e^{\eta_j\Phi_{k,j}}}. \tag{108}$$

Concerning the first term of (49), we have

$$\begin{aligned}
&\frac{\partial}{\partial\Phi_{a,j}}\left\{-\sum_{l=1}^{N_y}g\left(x_i,y_l,\frac{1}{\delta_l+1}\frac{\sum_{k=1}^{N_x}h\left(x_k\right)e^{\eta_l\Phi_{k,l}}\Delta x}{\sum_{k=1}^{N_x}e^{\eta_l\Phi_{k,l}}\Delta x}\right)W_{l,j}\Delta y\right\}\\
&=-\sum_{l=1}^{N_y}\frac{1}{\delta_l+1}g'\frac{\partial}{\partial\Phi_{a,j}}\left\{\frac{\sum_{k=1}^{N_x}h\left(x_k\right)e^{\eta_l\Phi_{k,l}}\Delta x}{\sum_{k=1}^{N_x}e^{\eta_l\Phi_{k,l}}\Delta x}\right\}W_{l,j}\Delta y \qquad (109)\\
&=-\frac{1}{\delta_j+1}g'_{l=j}\frac{\partial}{\partial\Phi_{a,j}}\left\{\frac{\sum_{k=1}^{N_x}h\left(x_k\right)e^{\eta_j\Phi_{k,j}}\Delta x}{\sum_{k=1}^{N_x}e^{\eta_j\Phi_{k,j}}\Delta x}\right\}W_{j,j}\Delta y,
\end{aligned}$$

where $g'$ denotes the partial derivative of $g$ with respect to its third argument (other arguments are omitted for simplicity). We have

$$\frac{\partial}{\partial\Phi_{a,j}}\left\{\frac{\sum_{k=1}^{N_x}h\left(x_k\right)e^{\eta_j\Phi_{k,j}}\Delta x}{\sum_{k=1}^{N_x}e^{\eta_j\Phi_{k,j}}\Delta x}\right\}=\frac{\eta_je^{\eta_j\Phi_{a,j}}}{\left(\sum_{k=1}^{N_x}e^{\eta_j\Phi_{k,j}}\right)^2}\left[h\left(x_a\right)\sum_{k=1}^{N_x}e^{\eta_j\Phi_{k,j}}-\sum_{k=1}^{N_x}h\left(x_k\right)e^{\eta_j\Phi_{k,j}}\right], \tag{110}$$

and hence

$$\begin{aligned}
&\frac{\partial}{\partial\Phi_{a,j}}\left\{-\sum_{l=1}^{N_y} g\left(x_i, y_l, \frac{1}{\delta_l+1}\frac{\sum_{k=1}^{N_x} h(x_k)e^{\eta_l\Phi_{k,l}}\Delta x}{\sum_{k=1}^{N_x} e^{\eta_l\Phi_{k,l}}\Delta x}\right)W_{l,j}\Delta y\right\} \\
&= -\frac{1}{\delta_j+1} g'_{l=j} \frac{\eta_j e^{\eta_j\Phi_{a,j}}}{\left(\sum_{k=1}^{N_x} e^{\eta_j\Phi_{k,j}}\right)^2}\left[h(x_a)\sum_{k=1}^{N_x} e^{\eta_j\Phi_{k,j}} - \sum_{k=1}^{N_x} h(x_k)e^{\eta_j\Phi_{k,j}}\right]W_{j,j}\Delta y \\
&\le \frac{1}{\delta_j+1} L_g \frac{\eta_j e^{\eta_j\Phi_{a,j}}}{\left(\sum_{k=1}^{N_x} e^{\eta_j\Phi_{k,j}}\right)^2}\left[\bar{h}\sum_{k=1}^{N_x} e^{\eta_j\Phi_{k,j}} - 0\right]W_{j,j}\Delta y \\
&= \frac{\bar{h}}{\delta_j+1} L_g \eta_j W_{j,j}\Delta y \frac{e^{\eta_j\Phi_{a,j}}}{\sum_{k=1}^{N_x} e^{\eta_j\Phi_{k,j}}}.
\end{aligned} \tag{111}$$

Consequently, we arrive at

$$\begin{aligned}
\frac{\partial\mathbb{S}}{\partial\Phi_{a,j}} &\le \frac{\bar{h}}{\delta_j+1} L_g \eta_j W_{j,j}\Delta y \frac{e^{\eta_j\Phi_{a,j}}}{\sum_{k=1}^{N_x} e^{\eta_j\Phi_{k,j}}} - \frac{1}{\delta_j}\frac{e^{\eta_j\Phi_{a,j}}}{\sum_{k=1}^{N_x} e^{\eta_j\Phi_{k,j}}} \\
&= \frac{e^{\eta_j\Phi_{a,j}}}{\sum_{k=1}^{N_x} e^{\eta_j\Phi_{k,j}}}\left(\frac{\bar{h}}{\delta_j+1} L_g \eta_j W_{j,j}\Delta y - \frac{1}{\delta_j}\right) \\
&\le 0,
\end{aligned} \tag{112}$$

the last line being due to (50).

Second, fix $(a,b)$ with $1 \le a \le N_x$ and $1 \le b \le N_y$ with $b \ne j$. For the second term of (49),

$$\frac{\partial}{\partial\Phi_{a,b}}\left\{-\frac{1}{\delta_j\eta_j}\ln\left(\sum_{k=1}^{N_x} e^{\eta_j\Phi_{k,j}}\Delta x\right)\right\} = 0 \tag{113}$$

because $b \ne j$. Concerning the first term of (49), as for (111), we have

$$\begin{aligned}
&\frac{\partial}{\partial\Phi_{a,b}}\left\{-\sum_{l=1}^{N_y} g\left(x_i, y_l, \frac{1}{\delta_l+1}\frac{\sum_{k=1}^{N_x} h(x_k)e^{\eta_l\Phi_{k,l}}\Delta x}{\sum_{k=1}^{N_x} e^{\eta_l\Phi_{k,l}}\Delta x}\right)W_{l,j}\Delta y\right\} \\
&= -\frac{1}{\delta_j+1} g'_{l=b} W_{b,j}\Delta y \frac{\eta_b e^{\eta_b\Phi_{a,b}}}{\left(\sum_{k=1}^{N_x} e^{\eta_b\Phi_{k,b}}\right)^2}\left[h(x_a)\sum_{k=1}^{N_x} e^{\eta_b\Phi_{k,b}} - \sum_{k=1}^{N_x} h(x_k)e^{\eta_b\Phi_{k,b}}\right],
\end{aligned} \tag{114}$$

and hence

$$\frac{\partial\mathbb{S}}{\partial\Phi_{a,b}} = -\frac{\eta_b}{\delta_b+1} W_{b,j}\Delta y g'_{l=b} \frac{e^{\eta_b\Phi_{a,b}}}{\sum_{k=1}^{N_x} e^{\eta_b\Phi_{k,b}}}\left[h(x_a) - \frac{\sum_{k=1}^{N_x} h(x_k)e^{\eta_b\Phi_{k,b}}}{\sum_{k=1}^{N_x} e^{\eta_b\Phi_{k,b}}}\right]. \tag{115}$$

The term inside the last “[ ]” does not necessarily have a fixed sign because both terms inside it may take any values between $\inf_{x\in\Omega} h(x)$ and $\sup_{x\in\Omega} h(x)$.

□

**A.2 Auxiliary computational results**

This appendix presents computational results for the solutions $\Phi$ to the HJB system (23). **Figure A1** shows the computed $\Phi$ for each case in **Figure 4**, while **Figure A2** shows the corresponding results for each case in **Figure 8**. The results imply that $\Phi$ remains continuous even when $\delta$ is less than 1, indicating that the assumptions in **Propositions 2-3** are sufficient but not necessary. Comparing **Figures A1 and A2**, we observe that the profiles of $\Phi$ vary more gradually in the no-graphon case than in the with-graphon case. This difference is attributed to the averaging effect of the graphon $W$ in the utility (32).

We finally report the convergence of numerical solutions to the HJB system for Case D in the with-Graphon case. **Table A1** reports the error of numerical solutions with $N_x = N_y = 2^M$ ($M = 4,5,6,7,8$) against a "reference" solution, where the numerical solution with $M = 9$ is considered the reference solution because no closed-form solution to the HJB system has been found. **Table A1** suggests the first- to second-order convergence of numerical solutions, which is consistent with the theoretical results presented in **Remark 3**. The iteration count does not significantly depend on the computational complexity $M$ and hence $N$. The computational time on a common laptop (Intel(R) Core(TM) Ultra 7 165H (3.80 GHz), 32.0 GB, Dev C++ Version 6.3) is regressed as $\log_2(\text{error}) = 2.58M - 13.93$ ($R^2$ value is 0.999), suggesting that the error grows exponentially in the grid size: $\text{error} = O\left(N^{2.58}\right)$. The computational complexity is larger than $O\left(N^2\right)$ that is expected for two-dimensional steady partial differential equations because ours is seen as a nonlocal equation.

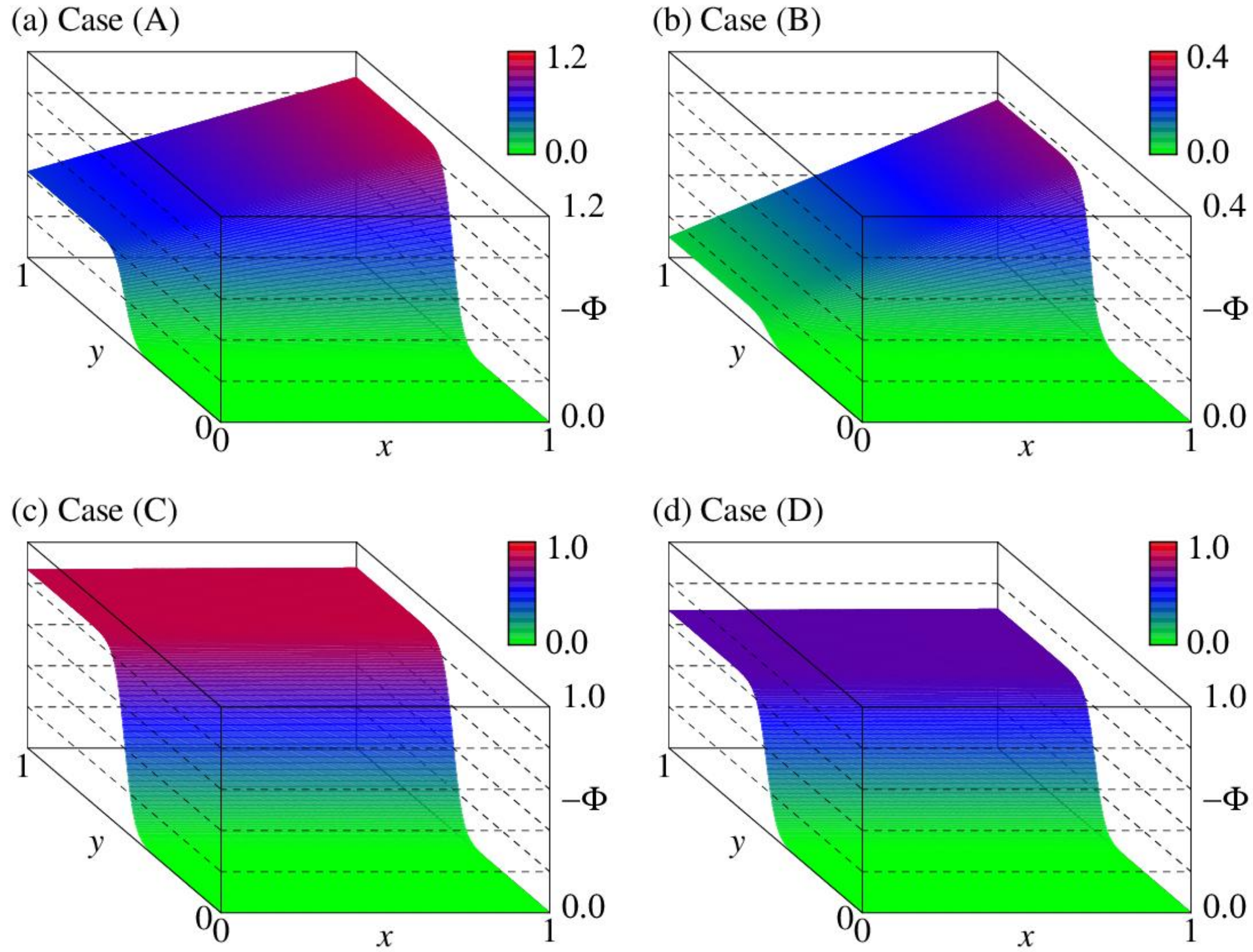

**Figure A1.** Computed $\Phi$ for Cases (A)-(D) for the G-MFLD (No-graphon case). We are plotting $-\Phi$ for visualization purposes: Case (A) $(\delta,\eta)=(0.5,1/0.5)$, Case (B) $(\delta,\eta)=(0.5,1/0.005)$, Case (C) $(\delta,\eta)=(0.005,1/0.5)$, and Case (D) $(\delta,\eta)=(0.005,1/0.005)$.

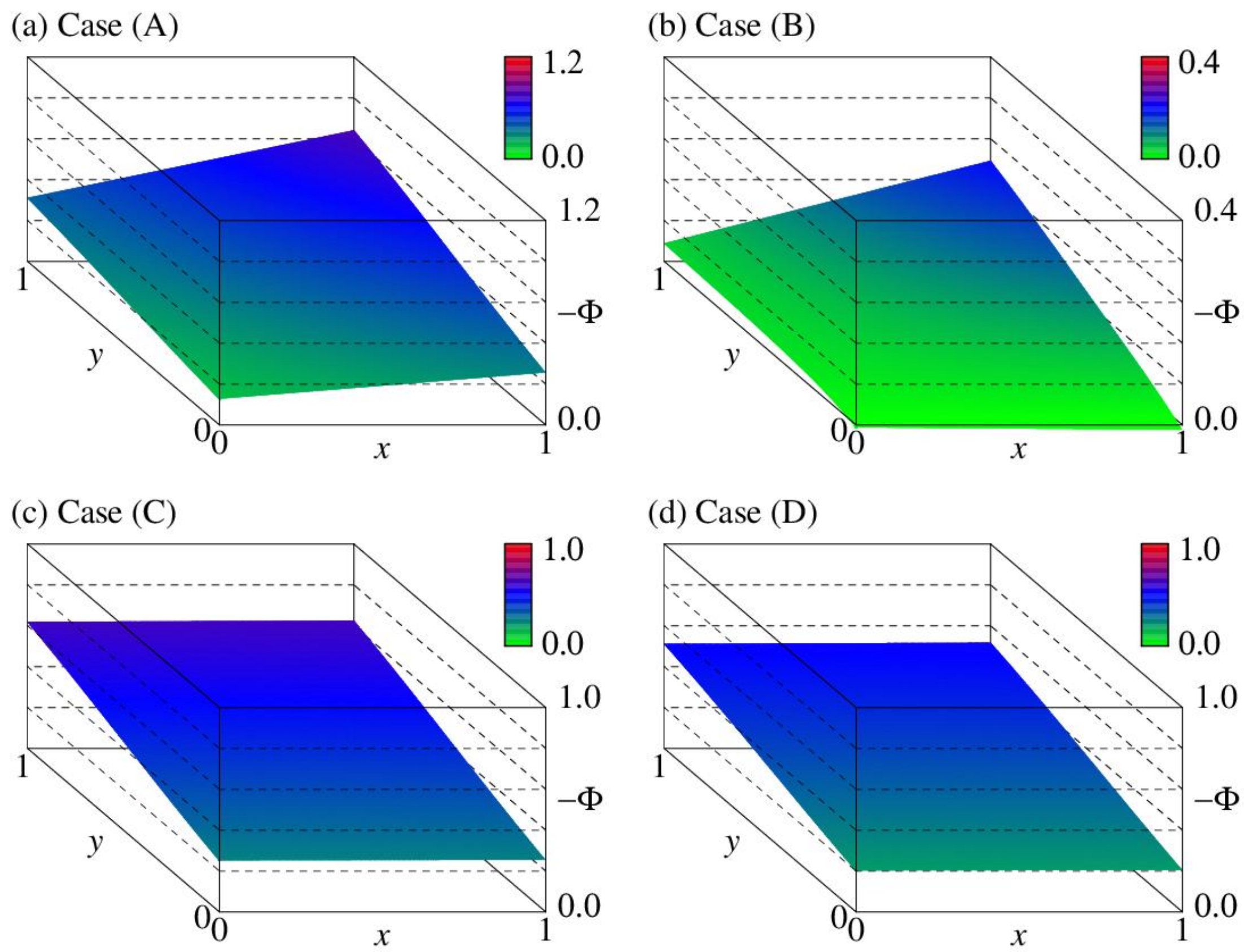

**Figure A2.** Computed solution $\Phi$ to the HJB system for Cases (A)-(D) for the G-MFLD (With-graphon case). We are plotting $-\Phi$ for visualization purposes: Case (A) $(\delta,\eta)=(0.5,1/0.5)$, Case (B) $(\delta,\eta)=(0.5,1/0.005)$, Case (C) $(\delta,\eta)=(0.005,1/0.5)$, and Case (D) $(\delta,\eta)=(0.005,1/0.005)$.

**Table A1.** Maximum pointwise error between numerical and reference solutions. Convergence rate at level $M$ is computed as $\log_2\left(\frac{\text{Error}|_{\text{Level } M}}{\text{Error}|_{\text{Level } M+1}}\right)$. We also report the iteration count and computational time (average of three runs) for each case.

| $M$ | Error | Convergence rate | Iteration count | Computational time (s) |
|---|---|---|---|---|
| 4 | 0.1510 | 1.153 | 2904 | 0.089 |
| 5 | 0.0679 | 1.294 | 2894 | 0.464 |
| 6 | 0.0277 | 1.704 | 2891 | 2.708 |
| 7 | 0.0085 | 2.087 | 2891 | 16.815 |
| 8 | 0.0020 | | 2891 | 112.260 |

**Declaration of interest** The author declares no relevant financial or non-financial interests.
**Declaration of Generative AI and AI-assisted technologies in the writing process** The author did not use generative AI in the writing of this manuscript.
**Data availability statement** The data will be made available upon reasonable request to the corresponding author.
**Funding** This study was supported by the Nippon Life Insurance Foundation [Environmental Problems Research Grant for Young Researchers, No. 24] and the Japan Society for the Promotion of Science [KAKENHI, No. 25K00240].
**Acknowledgments** The author would like to express his gratitude to the Hakusan Shiramine Fisheries Cooperative for providing valuable information on this study. The author gratefully thanks to helpful comments from the reviewers.